\documentclass[hidelinks,onefignum,onetabnum]{siamart250211}

\usepackage{verbatim}
\usepackage{latexsym, graphicx, epsfig}
\usepackage{pgfplots}
\usepackage{url}
\usepackage{mathrsfs}
\usepackage{color}
\usepackage{array}
\usepackage{amsfonts,amsmath,amssymb}
\usepackage{float}
\usepackage{multirow}
\usepackage{booktabs}
\usepackage{subcaption}
\usepackage{tikz}
\usepackage{diagbox}
\usepackage[toc,page]{appendix}
\usepackage{algorithm}
\usepackage{algpseudocode}
\usepackage[numbers,sort&compress]{natbib}
\usepackage{cleveref}

\newsiamremark{remark}{Remark}
\newsiamthm{assumption}{Assumption}
\numberwithin{equation}{section}
\numberwithin{figure}{section}

\title{An efficient Asymptotic-Preserving scheme for the Boltzmann mixture with disparate mass
\thanks{Submitted to the editors Nov 20, 2024.
\funding{Z. Hao acknowledges the computational resource of The Chinese University of Hong Kong during his visit.
N. Jiang acknowledges the support by NSFC grants 12371224, 11971360, 11731008 and the Strategic Priority Research Program of Chinese Academy of Sciences grant XDA25010404.
L. Liu acknowledges the support by National Key R\&D Program of China (2021YFA1001200), Ministry of Science and Technology in China, Early Career Scheme (24301021) and General Research Fund (14303022 \& 14301423) funded by Research Grants Council of Hong Kong.}}}

\author{
Zhen Hao \thanks{School of Mathematics and Statistics, Wuhan University, Wuhan 430072, P. R. China (\email{zhhao\_math@whu.edu.cn}).}
\and 
Ning Jiang \thanks{School of Mathematics and Statistics, Wuhan University, Wuhan 430072, P. R. China (\email{njiang@whu.edu.cn}).}
\and 
Liu Liu \thanks{The Chinese University of Hong Kong, Hong Kong 
  (\email{liuliu@math.cuhk.edu.hk}).}
}

\headers{AP Scheme for Disparate Mass Gas Mixtures}{Z. Hao, N. Jiang, and L. Liu}

\makeatletter
\providecommand{\headerps@out}{}
\makeatother

\begin{document}

\maketitle

\begin{abstract}
In this paper, we develop and implement an efficient asymptotic-preserving (AP) scheme to solve the gas mixture of Boltzmann equations under the disparate mass scaling relevant to the so-called ``epochal relaxation'' phenomenon. 
The disparity in molecular masses, ranging across several orders of magnitude, leads to significant challenges in both the evaluation of collision operators and the designing of time-stepping schemes to capture the multi-scale nature of the dynamics. 
A direct implementation of the spectral method faces prohibitive computational costs as the mass ratio increases due to the need to resolve vastly different thermal velocities. Unlike [I. M. Gamba, S. Jin, and L. Liu, \textit{Commun. Math. Sci.}, 17 (2019), pp. 1257-1289], we propose an alternative approach based on proper truncation of asymptotic expansions of the collision operators, which significantly reduces the computational complexity and works well for small $\varepsilon$. 
By incorporating the separation of three time scales in the model's relaxation process [P. Degond and B. Lucquin-Desreux, \textit{Math. Models Methods Appl. Sci.}, 6 (1996), pp. 405-436], we design an AP scheme that captures the specific dynamics of the disparate mass model while maintaining computational efficiency. Numerical experiments demonstrate the effectiveness of the proposed scheme in handling large mass ratios of heavy and light species, as well as capturing the epochal relaxation phenomenon.

\end{abstract}

\ifpdf
\hypersetup{
  pdftitle={AP scheme for a disparate mass gas mixture},
}
\fi

\begin{keywords}
gas mixture, Boltzmann equations, disparate mass, asymptotic analysis for collision operators, asymptotic-preserving scheme 
\end{keywords}

\begin{MSCcodes}
82C40, 65M70, 65T40, 65Y20, 35B25
\end{MSCcodes}

\section{Introduction}

The numerical computation of gas mixtures with disparate masses is a challenging task, with important applications in plasma physics and aerospace engineering. 
In practice, mass ratios of light and heavy molecules usually span several orders of magnitude, from  $10^2$ to $10^5$, as seen in spacecraft re-entry plasma environments \cite{Brun2012}, ITER fusion reactor dust-plasma interactions \cite{Takase}, and evaporation/condensation processes \cite{Takata-Golse}.

In kinetic theory, the Boltzmann equation models the evolution of rarefied gases and their mixtures. 
Major challenges for solving the Boltzmann equation include the non-local collision operator of an integral type and treatment of the multiple scales. 
Over the decades, deterministic methods, being free of statistical noise, have proven particularly advantageous for solving the Boltzmann equation, especially in near-continuum flows \cite{Dimarco-Pareschi}. In particular, the Fourier spectral method has been popularly used since the pioneering work \cite{Pareschi-Perthame, Pareschi-Russo}, with fast spectral methods developed afterwards, see, for example, \cite{Mouhot-Pareschi, Filbet-Mouhot-Pareschi, Gamba-Haack, Gamba-Haack-Hauck-Hu}. 
For gas mixtures, recent studies have achieved notable progress \cite{Jaiswal-Alexeenko-Hu, WuLei}.

Compared to single-species or standard multi-species cases, gas mixtures with disparate masses bring substantial computational challenges in the handling of collision operators. A key characteristic of these mixtures is the large disparity in thermal velocities of different species: 
the grid spacing must resolve both the fastest and the slowest species.
As a result, direct implementations of Fourier spectral methods
would have a computational complexity that scales with the mass ratio
\cite{Jaiswal-Alexeenko-Hu, WuLei}. Various studies address these challenges through adaptive meshing in velocity space \cite{Taitano-Chacon-Simakov} and asymptotic models tailored to large mass ratios \cite{Charles-Desvillettes, Taitano-Chacon-Simakov}.

Back in the 1960s, early efforts were made to derive the hydrodynamic models for plasmas and binary gas mixtures, the \textit{epochal relaxation} phenomenon in disparate mass mixtures was first pointed out by Grad \cite{Grad}, who remarked that the relaxation rate of the light species is faster than that of the heavy species by a factor depending on the square-root of the mass ratio $\varepsilon = \sqrt{m^L/m^H}$ where $L$ and $H$ denote the light and heavy species, respectively. 
There have been several analysis works along this line, see 
\cite{Braginskii,Spitzer-Harm,Goldman-Sirovich,Chmieleski-Ferziger,Johnson,Petit-Darrozes,Degond-Lucquin-Desreux}. In particular, Degond and Lucquin-Desreux \cite{Degond_bookchapter, Degond-Lucquin-Desreux, Degond-Lucquin-Desreux-2} characterized the epochal relaxation process by identifying three distinct time scales described by $\tau = 1$, $\tau = \varepsilon$ and $\tau = \varepsilon^2$. 
In contrast to normal mixture models, the disparate mass mixture exhibits scale separation at both kinetic and macroscopic levels, which brings significant computational challenges besides calculation for collision operators. Specifically, at the kinetic level, the distribution functions of different species {equilibrate} on distinct time scales; however at the macroscopic level, the velocities and temperatures of different species relax on separate time scales.

To capture the epochal relaxation phenomenon in numerical simulation, an efficient scheme not only needs to handle the stiffness brought by the small parameters $\varepsilon$ and $\tau$, but also allow efficient transition between kinetic and macroscopic regimes. 
To this end, we mention Asymptotic-Preserving (AP) schemes~\cite{Jin-AP}, whose goal is to avoid resolving the small scaling parameter numerically while preserving the asymptotic limits from the kinetic to the macroscopic models in the discrete setting. Such schemes allow \textit{automatic} switches from a kinetic solver to the macroscopic solver.
Several AP schemes have been designed for the Boltzmann mixture model, for examples see \cite{Jin-Shi,Jin-Li}. However, none of those work considered the epochal relaxation dynamics of disparate mass models.

The main focus of this work is to 
design an AP scheme for the Boltzmann mixture model with disparate masses that can capture the epochal relaxation phenomenon. 
To deal with the velocity disparity in the computation of inter-particle collision operators, we propose a novel approach that is based on proper truncation on the asymptotic expansion given in \cite{Degond_bookchapter,Degond-Lucquin-Desreux}. This reduces the computation to the evaluation of several velocity space integrals and derivatives. Similar approach has been used in~\cite{Taitano-Chacon-Simakov}.
To deal with the time stiffness, the BGK-penalty idea~\cite{Filbet-Jin} is employed. We mention that in \cite{Gamba-Jin-Liu}, the authors studied the same model and employed the ansatze $f^{L,H} = f_0^{L,H} + \varepsilon f_1^{L,H}$ to formulate a system of equations for the variables $f_0^{L,H}$ and $f_1^{L,H}$. 
They provided a theoretical AP proof of their method but lacked a practical implementation. Plus, their scheme still relies on the Fourier spectral method for the evaluation of inter-particle collision operators, which is computationally infeasible in practice.

\vspace*{0.5em}

\textbf{Novelty and contributions:} 
Compared to the previous work~\cite{Gamba-Jin-Liu} which lack of numerical simulations, we propose a significantly simpler approach and provide numerical experiments for Maxwell molecules. To the best of our knowledge, this is the {\it first} work in developing an efficient asymptotic-preserving scheme that can capture the so-called ``epochal relaxation" phenomenon for the Boltzmann mixture model with disparate mass, without suffering from the time stiffness when the square-root of the mass ratio $\varepsilon \ll 1$. In particular, we first derive the truncated model based on the asymptotic expansion in \cite{Degond_bookchapter,Degond-Lucquin-Desreux} and adopt the BGK-penalization technique by \cite{Filbet-Jin}, which requires updating the macroscopic quantities in the moment equations. Then, we carefully design an implicit scheme yet can be implemented explicitly for the moment equations, while ensuring stability and being able to capture the epochal relaxation phenomenon. Lastly, based on a careful analysis, we show the consistency of {the moment update} and the asymptotic-preserving property of our kinetic scheme.

Another major novelty and contribution of our work is developing a {\it uniformly} efficient numerical approach to evaluate the inter-particle collision operators when the mass ratio is small. 
Inspired by the expansion in \cite{Degond-Lucquin-Desreux} but from a numerical perspective, we design an effective way to compute the ``truncated" inter-particle collision operators and solve the mixture model. 
Specifically, we invent a dual-grid procedure to calculate the angular integrals for the asymptotic light-heavy collision operators, which involves a careful polar-grid design and interpolation between the polar and Cartesian grids. To this end, we provide an alternative and faster method compared to the existing expensive Fourier spectral scheme for disparate mass mixtures (see Remark 3 in~\cite{Jaiswal-Alexeenko-Hu}).

The rest of the paper is organized as follows: In Section \ref{sec:model}, we introduce the Boltzmann gas mixture with disparate mass scaling and list some key properties. Section \ref{sec:AE} describes the asymptotic expansion method for evaluating inter-particle collision operators, the polar grid design and angular integration, in addition to the applicability of our \textbf{AE} method in small $\varepsilon$ regime.
In Section \ref{sec:AP}, we present the time discretization, show the consistency of moments update and most importantly, the AP property of our numerical scheme. In Section \ref{sec:numerics}, we demonstrate the robustness and efficiency of the proposed scheme by several numerical experiments, which focus on the effectiveness of the asymptotic expansion and how the epochal relaxation phenomenon is captured. The conclusion is summarized in Section \ref{sec:conclusion}.

\section{The disparate-mass gas mixture model}
\label{sec:model}

In this section, we introduce the homogeneous Boltzmann equations for gas mixtures under the disparate mass scaling. We discuss the key properties of the inter-particle collision operators.
Finally, we review the ``epochal relaxation" phenomenon that is particular to disparate mass models described in \cite{Degond-Lucquin-Desreux}.

\subsection{The Boltzmann equation in disparate mass scaling}

The homogeneous Boltzmann equation for a binary gas mixture in the disparate mass scaling is given by~\cite{Degond-Lucquin-Desreux, Degond_bookchapter}
{\small
\begin{equation}
    \label{equation:main_equation}
    \begin{aligned}
        \tau \partial_t f^{L} &= Q^{LL}(f^{L}, f^{L}) + Q^{LH}_\varepsilon(f^{L}, f^{H}), \\
        \tau \partial_t f^{H} &= \varepsilon \left[Q^{HH}(f^{H}, f^{H}) + Q^{HL}_\varepsilon(f^{H}, f^{L})\right],
    \end{aligned}
\end{equation}
}
where $f^{L}$ and $f^{H}$ are the velocity distribution functions of the light and heavy species, respectively. 
Here, \(Q^{LL}\) and \(Q^{HH}\) denote collisions within the same species (hereafter referred to as `intra-particle' collisions),
and \(Q^{LH}_\varepsilon\) and \(Q^{HL}_\varepsilon\) denote collisions between different species (hereafter referred to as `inter-particle' collisions). The explicit expressions of the collision operators are given in Appendix~\ref{appendix:collision_operator_formulae}. 
The disparate mass regime concerns the case when the square-root of the mass ratio, defined as
{\small
\begin{equation}
    \varepsilon = \sqrt{\frac{m^L}{m^H}},
\end{equation}
}
is small.
Here, $m^L$ and $m^H$ are the masses of the light and heavy species, respectively. $\tau$ is the time scale of the problem.

For a velocity distribution function $f$, we denote $n$, $u$, and $T$ as its number density, mean velocity, and temperature 
{\small
\begin{equation}
    \label{equation:moments}
    n = \int_{\mathbb{R}^{d_{v}}} f dv, \quad u = \frac{1}{n} \int_{\mathbb{R}^{d_{v}}} f v dv, \quad T = \frac{1}{d_v n} \int_{\mathbb{R}^{d_{v}}} f |v-u|^2 dv.
\end{equation}
}
For a given set of macroscopic quantities $U = (n, u, T)^{\top}$, we define the Maxwellian distribution function as
{\small
\begin{equation}
    \label{equation:maxwellian}
    \mathcal{M}_{U}(v) = \mathcal{M}_{n, u, T}(v) = \frac{n}{(2\pi T)^{d_v / 2}} \exp \left( - \frac{|v-u|^2}{2 T} \right).
\end{equation}
}In this work, we refer to a Maxwellian distribution $M$ as \textit{associated to} a distribution function $f$ if 
$\label{equation:maxwellian_def}
M(v) = \mathcal{M}_{U}(v)$ where $U$ is computed from $f$ by~\eqref{equation:moments}.

\begin{remark}
    \label{remark:heavy_species_macros}
    In the disparate mass scaling derived in \cite{Degond-Lucquin-Desreux}, the heavy species velocity and distribution function are rescaled as \(
        \tilde{v}^{H} = \varepsilon v^{H}, 
        \quad  \tilde{f}^{H}(\tilde{v}^{H}) = \frac{1}{\varepsilon^{d_{v}}} f^{H}(v^{H}).\)
    Therefore, after a change of variables in~\eqref{equation:moments}, the density, mean velocity and temperature of the heavy species are given by $n^H$, $\varepsilon u^H$, and $T^H$, respectively, where $n^H$, $u^H$, and $T^H$ are computed from $f^H$ by~\eqref{equation:moments}.
\end{remark}

\subsection{Properties of the inter-particle collision operators}
The inter-particle collision operators $Q^{LH}_{\varepsilon}$ and $Q^{HL}_{\varepsilon}$ play a crucial role in the study of the epochal relaxation phenomenon.
We review some of the key properties given in \cite{Degond-Lucquin-Desreux,Degond-Lucquin-Desreux-2, Degond_bookchapter}.
\begin{proposition}[\textbf{Conservation properties of $Q^{LH}_{\varepsilon}$ and $Q^{HL}_{\varepsilon}$}]
    {\small
    $$
        \int_{\mathbb{R}^{d_{v}}} Q^{LH}_{\varepsilon} dv =
        \int_{\mathbb{R}^{d_{v}}} Q^{HL}_{\varepsilon} dv = 
            \int_{\mathbb{R}^{d_{v}}} Q^{LH}_{\varepsilon} v + Q^{HL}_{\varepsilon} v dv = 
            \int_{\mathbb{R}^{d_{v}}}Q^{LH}_{\varepsilon} |v|^2 + \varepsilon Q^{HL}_{\varepsilon} |v|^2 dv = 0.
    $$
    }
\end{proposition}
This corresponds to the conservation of mass, total momentum and total energy. Moreover, they can be expanded in terms of $\varepsilon$ as follows.
\begingroup
\begin{proposition}[\textbf{Asymptotic expansion of $Q^{LH}_\varepsilon$ and $Q^{HL}_\varepsilon$}]
\label{prop:asymptotic_expansions}
    Let $f^L (v)$ and $f^H (v)$ be sufficiently smooth functions. Then we have
    {\small
    \begin{equation}
        \label{equation:asymptotic_expansions}
        \begin{aligned}
            Q^{LH}_{\varepsilon} = \sqrt{1+\varepsilon^2} 
            \left( Q^{LH}_{0} + \varepsilon Q^{LH}_{1} 
            + \mathcal{O}(\varepsilon^2)\right),\\
            Q^{HL}_{\varepsilon} = \sqrt{1+\varepsilon^2} 
            \left( Q^{HL}_{0} + \varepsilon Q^{HL}_{1} 
            + \mathcal{O}(\varepsilon^2)\right).
        \end{aligned}
    \end{equation}
    }These asymptotic operators $\{Q^{LH}_i\}_{i\in\mathbb{N}}$ and $\{Q^{HL}_i\}_{i\in\mathbb{N}}$ have the following properties: 
      \begin{enumerate}
        \item For any function $f^{H}$, we have
        \begin{enumerate}
            \item If $f^{L}$ is a function of $|v^{L}|$, then $Q^{LH}_{0}(f^{L},f^{H})=0$, 
            \item If $f^{L}$ is an even function, then $Q^{HL}_{0}(f^{H},f^{L})=0$.
        \end{enumerate}
        \item Conservation properties:
        {\small
        \begin{equation}
            \label{equation:conservation_properties}
            \begin{aligned}
                &\int_{\mathbb{R}^{d_{v}}}Q^{LH}_{i}\mathrm{d}v = \int_{\mathbb{R}^{d_{v}}}Q^{HL}_{i}\mathrm{d}v = 0,\\
                &\int_{\mathbb{R}^{d_{v}}}Q^{LH}_{i}v + Q^{HL}_{i}v\mathrm{d}v = 0,\quad \text{for } i\geqslant 0,\\
                &\int_{\mathbb{R}^{d_{v}}}Q^{LH}_{0}|v|^2\mathrm{d}v = 0, \quad \int_{\mathbb{R}^{d_{v}}}Q^{LH}_{i}|v|^2 + Q^{HL}_{i-1}|v|^2\mathrm{d}v = 0,\quad \text{for } i\geqslant 1.
            \end{aligned}
        \end{equation}
        }
        \end{enumerate}
\end{proposition}
\endgroup
We will derive the explicit formulae of the asymptotic operators in section~\ref{sec:AE}.

\subsection{Epochal relaxation and the macroscopic model}
\label{subsection:epochal_relaxation}
In a series of work by Degond and Lucquin-Desreux \cite{Degond-Lucquin-Desreux,Degond-Lucquin-Desreux-2}, 
an asymptotic analysis based on expansions of the inter-particle collision operators provides a concrete picture of the epochal relaxation phenomenon. Specifically, there is a three-time scale separation in the disparate mass regime. 
\begin{enumerate}
    \item[(i)] The fastest time scale ($\tau = 1$): collision time of the light species. 
    The heavy particle distribution does not evolve in time and the light particle one is described by a kinetic equation with two collision terms corresponding to self collisions and an elastic scattering of the light particles against the heavy ones as if the latter were steady.
    Also, the relaxation of velocities occurs at this time scale.
    \item[(ii)] The intermediate time scale ($\tau = \varepsilon$): collision time of the heavy species. 
    The light particles are at thermodynamic equilibrium with a centered Maxwellian distribution, while the heavy ones are subject to collisions with particles of the same species.
    \item[(iii)] The slowest time scale ($\tau = \varepsilon^2$): relaxation time scale. Both distributions are thermalized and the temperatures evolve one to each other via a relaxation equation.
    There exist $n^L, n^H \geq 0$, 
    $u^H \in \mathbb{R}^{d_v}$ and $T^L, T^H \geq 0$ such that
    {\small
    \begin{equation*}
            f^L_\varepsilon = \mathcal{M}_{n^L, 0, T^L} + \mathcal{O}(\varepsilon),\quad f^H_\varepsilon = \mathcal{M}_{n^H, u^H, T^H} + \mathcal{O}(\varepsilon). 
    \end{equation*}
    }Here the temperatures $(T^L(t), T^H(t))$ satisfy the following relaxation equations: 
    {\small
    \begin{equation}
            \label{equation:macro}
     \left\{
            \begin{aligned}
                    \frac{\mathrm{d}}{\mathrm{d}t}\left(\frac{d_v n^L T^L}{2}\right)
                    = -d_v \frac{\lambda (T^L)}{T^L}  n^H (T^L - T^H), \\[6pt]
                    \frac{\mathrm{d}}{\mathrm{d}t}\left(\frac{d_v n^H T^H}{2}\right)
                    = -d_v \frac{\lambda (T^L)}{T^L}  n^H (T^H - T^L),
      \end{aligned}\right.
    \end{equation}
    }where $\lambda (T)$ is given by
    {\small
    \begin{equation*}
        \lambda (T) = \frac{2}{d_v} \int_{\mathbb{R}^{d_v}} \int_{\mathbb{S}^{d_v-1}} B^{LH}(|v|, \sigma) |v|^2 {\mathcal{M}_{n^L, 0, T}} \mathrm{d}\sigma \mathrm{d}v
    \end{equation*}
    }
\end{enumerate}

\section{Efficient approximation of collision operators}
\label{sec:AE}

Numerical implementation of the Boltzmann collision operator has been a challenging problem due to the high dimension of its collision integral. 
To this end, a class of methods based on Fourier spectral approximations has been proposed in several works~\cite{Pareschi-Russo, Mouhot-Pareschi,Gamba-Haack-Hauck-Hu}. Among others, the Mouhot-Pareschi fast spectral method~\cite{Mouhot-Pareschi} is able to achieve a $\mathcal{O}(N^{d_v}\mathrm{log}N)$ time complexity, where $N$ denotes the number of points used in the velocity discretization. The intra-particle collision operators, being the same as the single-species Boltzmann collision operator (see~\eqref{equation:intra_collision_operators}), can be evaluated by the fast spectral method.

For the inter-particle collision operators $Q^{LH}_\varepsilon$ and $Q^{HL}_\varepsilon$, we mention the fast spectral methods by \cite{Jaiswal-Alexeenko-Hu, WuLei} that achieve a computational complexity of
$\mathcal{O}(N^{d_v+1}\mathrm{log}N)$.
However, as remarked in \cite{Jaiswal-Alexeenko-Hu, WuLei}, for large mass ratios, the spectral method is complicated by the non-unitary mass ratio between different molecular particles. Direct numerical implementation in disparate mass regimes requires a grid resolution increasing in the order of $\mathcal{O}(1/\varepsilon^{2d_v})$, resulting in a prohibitively high computational cost.

On the other hand, the scaling of the inter-particle collision operators here is different from those in \cite{Jaiswal-Alexeenko-Hu, WuLei}. 
We present in Appendix~\ref{appendix:SP} a modified  spectral method (hereafter refered to as the \textbf{SP} method) for $Q^{LH}_\varepsilon$ and $Q^{HL}_\varepsilon$ based on the work by Jaiswal-Alexeenko-Hu \cite{Jaiswal-Alexeenko-Hu}.

Inspired by the asymptotic expansions in section~\ref{sec:model}, we propose to truncate the asymptotic expansions of $Q^{LH}_\varepsilon$ and $Q^{HL}_\varepsilon$ to a certain order, and compute the truncated asymptotic operators as efficient approximations for the original operators. Since our primary focus is on the disparate mass regime where $\varepsilon \ll 1$, we expect the truncation to be accurate.
We remark that this approach resonates with a possible alternative to the Fourier spectral method in disparate mass regimes mentioned in~\cite{Jaiswal-Alexeenko-Hu}. In the remaining part of this section, we will present the details of the asymptotic-expansion (AE) method.

\subsection{The asymptotic-expansion (AE) method for $Q^{LH}_\varepsilon$ and $Q^{HL}_\varepsilon$}
\label{subsec:AE}

Recall from Proposition~\ref{prop:asymptotic_expansions} that one can derive the asymptotic expansions for $Q^{LH}_{\varepsilon}$ and $Q^{HL}_{\varepsilon}$: 
{\small
\begin{align*}
    Q^{LH}_{\varepsilon} =& \sqrt{1+\varepsilon^2} 
    \left( Q^{LH}_{0} + \varepsilon Q^{LH}_{1} 
    + \varepsilon^2 Q^{LH}_2 + \mathcal{O}(\varepsilon^3)\right),\\
    Q^{HL}_{\varepsilon} =& \sqrt{1+\varepsilon^2} 
    \left( Q^{HL}_{0} + \varepsilon Q^{HL}_{1} 
    + \mathcal{O}(\varepsilon^2)\right).
\end{align*}
}The AE method consists of approximating $Q^{LH}_\varepsilon$ and $Q^{HL}_\varepsilon$ by 
{\small
\begin{equation}
    \label{equation:AE}
    \begin{aligned}
        Q^{LH}_\varepsilon &\approx Q^{LH}_{\text{AE}} = \sqrt{1+\varepsilon^2}(Q^{LH}_0 + \varepsilon Q^{LH}_1 + \varepsilon^2 Q^{LH}_2),\\
        Q^{HL}_\varepsilon &\approx Q^{HL}_{\text{AE}} = \sqrt{1+\varepsilon^2}(Q^{HL}_0 + \varepsilon Q^{HL}_1),
    \end{aligned}
\end{equation}
}In the case $d_v=2$ and Maxwell molecules where $B^{LH}$, $B^{HL}$ are constant, the asymptotic collision operators can be written in the following simplified forms. The sketch of a derivation can be found in Appendix~\ref{appendix:derivation}. 
The light-heavy collision operators for $v^L\neq 0$ are given by
{\small
\begin{equation}
    \label{equation:QLH_012}
    \begin{aligned}
        Q^{LH}_0(v^L) &= B^{LH} n^H \left( \langle f^L \rangle - 2\pi f^L\right),\\
        \vspace*{3pt}
        Q^{LH}_1(v^L) &= B^{LH} n^H u^H \cdot 
        \left( \langle \nabla_{v^L} f^L \rangle - \nabla_{v^L} \langle f^L \rangle \right), \\
        Q^{LH}_2 (v^L) &= 
        B^{LH} \Big\{
            2n^H \langle f^L \rangle - n^H \frac{v^L}{|v^L|} \cdot \langle \sigma f^L \rangle 
            + n^H |v^L| \langle \sigma \cdot \nabla_{v^L} f^L \rangle \\
            & + \Big(\frac{1}{2} n^H |u^H|^2 + n^H T^H \Big)\frac{1}{|v^L|} \langle \sigma \cdot \nabla_{v^L} f^L \rangle 
            - n^H v^L \cdot \langle \sigma \otimes \sigma \cdot \nabla_{v^L} f^L \rangle \\
            & + \frac{1}{2} \int_{\mathbb{R}^2} v^H \otimes v^H f^H \mathrm{d} v^H : \Big(
                -\frac{v^L \otimes v^L}{|v^L|^3} \langle \sigma \cdot \nabla_{v^L} f^L \rangle 
                + \langle \nabla_{v^L}^2 f^L \rangle \\
                & - 2\langle \nabla_{v^L}^2 f^L \cdot \sigma \rangle \otimes \frac{v^L}{|v^L|}
                + \frac{v^L \otimes v^L}{|v^L|^2} \langle \sigma \otimes \sigma : \nabla_{v^L}^2 f^L \rangle
            \Big)\Big\}, 
    \end{aligned}
\end{equation}
}where $\langle f \rangle = \int_{\mathbb{S}^1} f(\sigma) \mathrm{d} \sigma$, and
$$ Q^{LH}_0 (0) = 0, \quad Q^{LH}_1 (0) = 2\pi B^{LH} n^H u^H \cdot \nabla_{v^L} f^L (0), \quad Q^{LH}_2 (0) = 4\pi B^{LH} n^H f^L (0). $$
The heavy-light collision operators are given by 
{\small
\begin{equation}
    \label{equation:QHL_01}
    \begin{aligned}
        Q^{HL}_0(v^H) =& -2\pi B^{HL}\,\nabla_{v^H}f^H \cdot n^L u^L,\\
        Q^{HL}_1(v^H) =& 2\pi B^{HL}\left( v^H \cdot \nabla_{v^H} f^H n^L + 2 f^H n^L \right)\\
        +& \pi B^{HL} \Delta_{v^H} f^H \left( \frac{1}{2} n^L |u^L|^2 + n^L T^L \right) \\
        +& \pi B^{HL} \nabla_{v^H}^2 f^H : \int_{\mathbb{R}^2} v^L \otimes v^L f^L \mathrm{d}v^L.
    \end{aligned}
\end{equation}
}
\begin{remark}
    The necessity to include the higher order term $Q^{LH}_{2}$ in the \textbf{AE} method lie in two aspects: (1) to attain the consistency of~\eqref{equation:main_equation} at the slowest time scale ($\tau = \varepsilon^2$); (2) from the analysis aspect, the dynamics of light-heavy species collisions at the slowest time scale is dominated by $Q^{LH}_{1}$ and $Q^{LH}_{2}$, as shown in \cite[equation (5.52)]{Degond-Lucquin-Desreux}.
\end{remark}

\subsection{Evaluation of collision operators in the AE method}

In this subsection, we will carefully study how to numerically implement the collision operators \eqref{equation:QLH_012} and \eqref{equation:QHL_01} appearing in the AE method. 
We consider two-dimensional velocity variable with computational domain $[-L_v, L_v]^2$. Assume that 
$f^L(v)$ and $f^H(v)$ are periodic in $v$, and the length $L_v$ is chosen within which $f^L(v)$ and $f^H(v)$ are compactly supported. The supports and $L_v$ satisfy the de-aliasing condition in \cite{Mouhot-Pareschi}. 

To compute $Q_0^{HL}$, $Q_1^{HL}$ in ~\eqref{equation:QHL_01}, one only needs simple operations of differentiation and integration, which can be computed by the central difference and trapezoidal rule. As for $Q_0^{LH}$, $Q_1^{LH}$, and $Q_2^{LH}$ in \eqref{equation:QLH_012}, the main challenge lies in the evaluation of angular integration. Since the Fourier spectral method for the intra-particle collision operators $Q^{LL}$ and $Q^{HH}$ requires a Cartesian grid, we will first interpolate the values of $f^L(v)$ into values $f^L_{\text{pol}}(r,\sigma)$ defined on a polar grid, then evaluate the light-heavy collision operators in polar coordinates, finally interpolate the values of collision operators back onto the Cartesian grid. The polar grid and interpolation of function values between the two grids will be discussed in the following subsection. 

\subsubsection{Polar grid and interpolation}

\begin{figure}[h]
    \centering
    \includegraphics[width=0.5\textwidth]{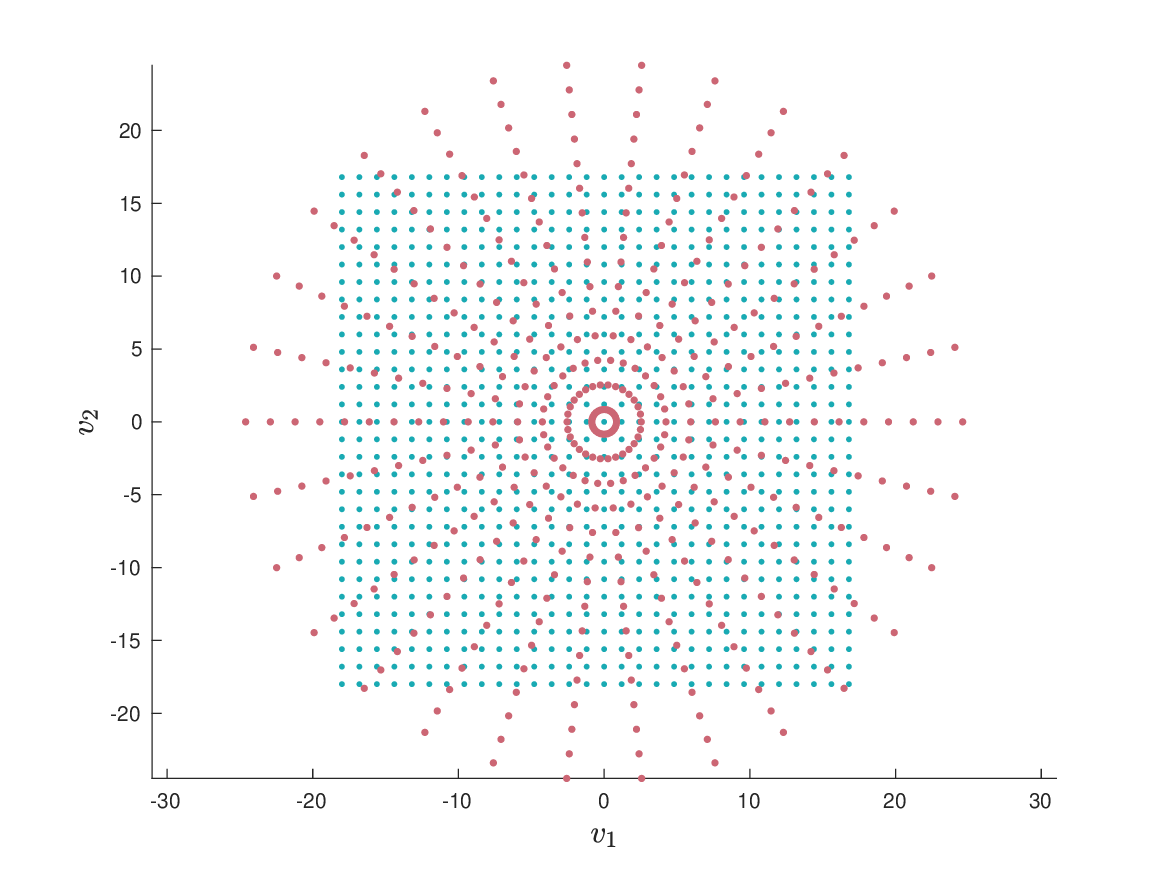}
    \caption{Illustration of grid design with $N_v = 30$}
    \label{figure:grid_design}
\end{figure}

For a non-zero vector \( v \), let \( v = r \sigma = r(\cos \theta, \sin \theta) \), where \( r \in (0, |v|_{\text{max}}] \) and \( \theta \in [0, 2\pi) \). The angle \( \theta \) is discretized periodically into $N_\theta$ grid points
\(
0 = \theta_1 < \theta_2 < \ldots < \theta_{N_\theta},
\)
where \( \theta_i = (i-1) \frac{2\pi}{N_\theta} \) for \( i = 1, \ldots, N_\theta \).
To avoid singularity at the origin, the radial component \( r \) is discretized into $N_r$ grid points defined by 
\(
\frac{\Delta r}{2} = r_1 < r_2 < \ldots < r_{N_r} = |v|_{\text{max}} - \frac{\Delta r}{2}, 
\)
where \( r_{j+1} - r_j = \Delta r \) for \( j = 1, \ldots, N_r - 1 \), and the mesh size is 
$\Delta r = \frac{|v|_{\text{max}}}{N_r}.$ 
We let $N_r = N_v/2$, $N_\theta = N_v$ and show the layout of a Cartesian grid and a polar grid with $N_v = 30$ in Figure~\ref{figure:grid_design}. 
In our simulations, much finer grids are used. 
We give some remarks about the Cartesian-Polar grid design. To begin with, the velocity domains covered by the two types of grids are different. In addition, the Cartesian grid is built on a uniform mesh, while the polar grid exhibit a clustering phenomenon near the origin and becomes sparser as the radius increases. As mentioned in \cite{Mouhot-Pareschi}, the Fourier spectral method for the intra-particle collision operators requires $S \approx 0.38 L_v$, with $S$ being the support of distribution functions. This indicates that the supports of $f^L$, $f^H$ are concentrated near the origin, thus the accuracy loss brought by the polar grid layout in the far-away region is relatively negligible, and it is safe to set the values on those ``suburb" points to be zero. 

\subsubsection{Collision operators in polar coordinates}
To compute the light-heavy collision operators in polar coordinates, we need to rewrite~\eqref{equation:QLH_012}. By a change of variable $v = r\sigma = r(\mathrm{cos}\theta, \mathrm{sin}\theta)$, $r\neq 0$, we first transform the derivatives in Cartesian coordinates into polar coordinates,
then apply integration by parts to eliminate derivatives in $\theta$ to get 
\begin{equation}
    \label{equation:QLH_012_polar}
    \begin{aligned}
    Q^{LH}_0(v^L) =& B^{LH} n^H \left[ \langle f^L \rangle - 2\pi f^L\right], \\
    Q^{LH}_1(v^L) =& B^{LH} n^H 
    \left\{ u^H_1 \left[  
            \langle \mathrm{cos}\theta\partial_r f^L \rangle 
            + \frac{1}{r} \langle \mathrm{cos}\theta f^L \rangle
            - \mathrm{cos}\theta \langle \partial_r f^L \rangle
        \right]\right. \\
        &\left. + u^H_2 \left[
            \langle \mathrm{sin}\theta\partial_r f^L \rangle 
            - \frac{1}{r} \langle \mathrm{sin}\theta f^L \rangle
            - \mathrm{sin}\theta \langle \partial_r f^L \rangle
        \right]
    \right\}, \\
    Q^{LH}_2(v^L) =& B^{LH} \{I_1 + I_2 + I_3 + I_4 + I_5\},
\end{aligned}
\end{equation}
where we define $ \langle f \rangle := \int_{\mathbb{S}^1} f(r\mathrm{cos}\theta, r\mathrm{sin}\theta) \mathrm{d} \theta $. With the following definitions
{\small
\begin{align*}
    A_1 &= \langle \partial_{rr} f^L \rangle - \frac{1}{r} \langle \partial_r f^L \rangle, \quad A_2 = \frac{1}{r} \langle \partial_r f^L \rangle, \quad A_5 = \langle \partial_r f^L \rangle, \\
    A_3 &= \partial_{rr} f^L + \frac{3}{r} \partial_r f^L + \frac{2}{r^2} f^L, \quad A_4 = \partial_{rr} f^L + \frac{1}{r} \partial_r f^L - \frac{1}{r^2} f^L,
\end{align*}
}
the components \( I_1 \) through \( I_5 \) are given by

{\small
\vspace{-10pt}
\begin{align*}
    I_1 &= 2n^H \langle f^L \rangle 
    - \cos\theta \langle \cos\theta f^L \rangle
    - \sin\theta \langle \sin\theta f^L \rangle, \; \\
% \end{align*}
% \vspace{-10pt}
% \begin{align*}
    I_2 &= n^H r A_5 
    + \left( \frac{1}{2} n^H |u^H|^2 + n^H T^H \right) A_2, \\
% \end{align*}
% \vspace{-10pt}
% \begin{align*}
    I_3 &= - n^H r \Big( \cos\theta \langle \cos\theta \partial_r f^L \rangle 
    + \sin\theta \langle \sin\theta \partial_r f^L \rangle \Big), \\
% \end{align*}
% \vspace{-10pt}
% \begin{align*}
    I_4 &= \frac{1}{2} \int (v^H_1)^2 f^H \, \mathrm{d}v^H \Big[ 
        \cos^2\theta A_1 
        + \langle \cos^2\theta \partial_{rr} f^L \rangle \\
        &\qquad + \frac{1}{r} \langle (2\cos^2\theta - \sin^2\theta) \partial_r f^L \rangle - 2\cos\theta \langle \cos\theta A_4 \rangle \Big], \\
% \end{align*}
% \vspace{-10pt}
% \begin{align*}
    I_5 &= \int v^H_1 v^H_2 f^H \, \mathrm{d}v^H \Big[ 
        \cos\theta \sin\theta A_1
        + \langle \cos\theta \sin\theta A_3 \rangle  - \cos\theta \langle \cos\theta A_4 \rangle \\
        &\qquad - \sin\theta \langle \sin\theta A_4 \rangle \Big] + \frac{1}{2} \int (v^H_2)^2 f^H \, \mathrm{d}v^H \Big[ 
        \sin^2\theta A_1
        + \langle \sin^2\theta \partial_{rr} f^L \rangle  \\
        &\qquad + \frac{1}{r} \langle (2\sin^2\theta - \cos^2\theta) \partial_r f^L \rangle - 2\sin\theta \langle \sin\theta A_4 \rangle \Big],
\end{align*}
}
Note that the numerical computations are reduced to integrals in $\theta$ and derivatives in $r$. 
Thanks to the polar-grid design, for the asymptotic operators shown in \eqref{equation:QLH_012_polar}, we can efficiently evaluate the angular integration by using the trapezoidal rule; while differentiation in $r$ can be computed by the standard central difference method. 

\subsection{Applicability of the \textbf{AE} method for small $\varepsilon$ regimes}
\label{subsec:AE_applicability}
Our designed \textbf{AE} method is based on the truncation of $Q^{LH}_{\varepsilon}$ and $Q^{HL}_{\varepsilon}$, introducing a local truncation error of $\mathcal{O}(\varepsilon^3/\tau)$ to the equation~\eqref{equation:main_equation}, which is small under all three time scales when $\varepsilon \ll 1$. Therefore, \textbf{AE} method is accurate only for small $\varepsilon$ regimes--this is fine since the main objective of this work is to study the dynamics of {\it disparate mass} mixtures. On the other hand, regarding the computational complexity, the \textbf{AE} method {requires $\mathcal{O}(N_v^{d})$ operations independently of $\varepsilon$}, while the traditional \textbf{SP} method on the non-truncated model requires $\mathcal{O}((N_v(\varepsilon))^{2d})$ as shown in Proposition~\ref{Prop:SP}. Figure~\ref{figure:AE_SP_regimes} illustrates the regimes where our \textbf{AE} method is more suitable and advantageous in terms of accuracy and efficiency. We will make further comparison between the two methods in the first numerical example below. 
\begin{figure}
    \centering
    \resizebox{0.7\textwidth}{!}{
    \begin{tikzpicture}[scale=1]

    \draw[<-, thick] (0, 0) -- (10, 0);
    
    \draw[thick] (0.65, 0.05) -- (0.65, -0.05) node[below]{$0.01$}; 
    \draw[thick] (1.45, 0.05) -- (1.45, -0.05) node[below]{$0.03$};
    \draw[thick] (6.5, 0.05) -- (6.5, -0.05) node[below]{$0.1$};   
    \draw[thick] (2.3, 0.05) -- (2.3, -0.05);
    \draw[thick] (3.3, 0.05) -- (3.3, -0.05);
    \draw[thick] (4.3, 0.05) -- (4.3, -0.05);
    \draw[thick] (8.25, 0.05) -- (8.25, -0.05) node[below]{$0.5$};       
    \draw[thick] (10, 0.05) -- (10, -0.05) node[below]{$1$};         

    % AE Method centered at x=2.5
    \node at (2.5, 0.9) {\small \textbf{AE} Method}; % Label slightly above arrows
    \draw[<-, thick] (0.6, 0.6) -- (1.6, 0.6) % Arrow pointing left
        node[midway, below, sloped]{\footnotesize efficient};
    \draw[->, thick] (3, 0.6) -- (4, 0.6) % Arrow pointing right
        node[midway, below, sloped]{\footnotesize lose accuracy};

    % SP Method centered at x=7.5
    \node at (7.5, 0.9) {\small \textbf{SP} Method}; % Label slightly above arrows
    \draw[<-, thick] (6, 0.6) -- (7, 0.6) % Arrow pointing left
        node[midway, below, sloped]{\footnotesize expensive};
    \draw[->, thick] (8, 0.6) -- (9, 0.6) % Arrow pointing right
        node[midway, below, sloped]{\footnotesize accurate};

    \node at (5, -0.8) {$\varepsilon$ values};

\end{tikzpicture}
    }
\caption{\footnotesize Regimes where the \textbf{AE} method is best suited. In large $\varepsilon$ regimes, the \textbf{SP} method is accurate, but the \textbf{AE} method loses accuracy; in small $\varepsilon$ regimes, the \textbf{SP} method requires prohibitively fine grid resolution while the \textbf{AE} method is accurate and efficient.}
\label{figure:AE_SP_regimes}
\end{figure}

\section{Asymptotic-preserving (AP) time discretization}
\label{sec:AP}

In this section, we introduce the asymptotic-preserving (AP) time discretization of \eqref{equation:main_equation} with the \textbf{AE} method used to approximate the inter-particle collision operators.
The AP scheme is based on the BGK-penalization technique first proposed by Filbet and Jin \cite{Filbet-Jin}, and inspired by the study in \cite{Gamba-Jin-Liu}. 
We first penalize the collision operators by the light species and heavy species BGK operators,

{\footnotesize
\begin{equation}
    \label{equation:scheme_initial}
    \begin{aligned}
        \frac{f^{n+1}_{L}-f^{n}_{L}}{\Delta t}
        =& \frac{1}{\tau}\left(Q^{LL}(f^n_L,f^n_L) + Q^{LH}_{\text{AE}}(f^n_L,f^n_H)
        - \nu^{n}_{L}(M^n_{L}-f^{n}_{L})\right) \\ & + \frac{1}{\tau}\nu^{n+1}_{L}(M^{n+1}_{L}-f^{n+1}_{L}), \\
        \frac{f^{n+1}_{H}-f^{n}_{H}}{\Delta t}
        =& \frac{\varepsilon}{\tau}\left(Q^{HH}(f^{n}_{H},f^{n}_{H}) + Q^{HL}_{\text{AE}}(f^n_H,f^n_L)
        -\nu^{n}_{H}(M^n_{H}-f^{n}_{H})\right) \\ & + \frac{\varepsilon}{\tau} \nu^{n+1}_{H} (M^{n+1}_{H}-f^{n+1}_{H}),
    \end{aligned}
\end{equation}
}where $n$ stands for the time step and the computational time $t = n \Delta t$. Here $M^{n, n+1}_L$ and $M^{n, n+1}_H$ are Maxwellians associated to $f^{n, n+1}_L$ and $f^{n, n+1}_H$ according to~\eqref{equation:maxwellian_def}.

One of the advantages of the BGK-penalization approach lies in the possibility of obtaining $M^{n+1}_L$ and $M^{n+1}_H$ \textit{a priori}, making it possible to treat the implicit method \eqref{equation:scheme} \textit{explicitly}. 
The typical approach to obtain $M^{n+1}_L$ and $M^{n+1}_H$ is to integrate~\eqref{equation:scheme} in $v$ against $\phi(v) = 1,\, v,\, |v|^2$, resulting in a system of equations for the macroscopic unknowns. The solution to the system is then used to approximate the moments of $f^{n+1}_L,\, f^{n+1}_H$ (see definition~\eqref{equation:moments}), which in turn are used to define $M^{n+1}_L$ and $M^{n+1}_H$.
In the model under consideration, a direct integration of~\eqref{equation:scheme} gives
{\small
\begin{equation}
    \label{equation:taking_moments}
    \begin{aligned}
        & \frac{U^{n+1}_L - U^n_L}{\Delta t}  =  \langle \frac{f^{n+1}_{L}-f^{n}_{L}}{\Delta t}, \phi(v^L) \rangle \\
        &\quad = \frac{1}{\tau} \langle Q^{LH}_0(f^n_L,f^n_H), \phi(v^L) \rangle + \frac{\varepsilon}{\tau} \langle Q^{LH}_1(f^n_L ,f^n_H), \phi(v^L) \rangle + \frac{\varepsilon^2}{\tau} \langle Q^{LH}_2(f^n_L,f^n_H), \phi(v^L) \rangle \\
        &\frac{U^{n+1}_H - U^n_H}{\Delta t} =  \langle \frac{f^{n+1}_{H}-f^{n}_{H}}{\Delta t}, \phi(v^H) \rangle \\ &\quad = \frac{\varepsilon}{\tau}\langle Q^{HL}_{0}(f^n_H, f^n_L), \phi(v^H) \rangle + \frac{\varepsilon^2}{\tau}\langle Q^{HL}_{1}(f^n_H,f^n_L), \phi(v^H) \rangle,
    \end{aligned}
\end{equation}
}where $\langle f, \phi \rangle := \int_{\mathbb{R}^{d_v}} f(v) \phi(v) \mathrm{d}v$,
and $U^n_L,\, U^n_H$ are the moments of $f^n_L,\, f^n_H$. 
The moment system is still stiff when $\tau \leq \varepsilon$, an issue also observed in \cite{GJL19}. However, for Maxwell molecules, explicit expressions can be computed from~\eqref{equation:QLH_012} and~\eqref{equation:QHL_01}. This is stated in the following lemma.
\begin{lemma}
    \label{lemma:collision_explicit_formula}    
    For Maxwell molecules, we have the following relations
    {\small
    \begin{align*}
        &\langle Q^{HL}_0(f^H, f^L), v^H \rangle = 2 \pi B^{HL} n^H n^L u^L, \; \langle Q^{HL}_0(f^H, f^L), |v^H|^2 \rangle = 4 \pi B^{HL} n^H n^L u^L \cdot u^H, \\
        &\langle Q^{HL}_1(f^H, f^{L}), v^H \rangle = -2 \pi B^{HL} n^H n^L u^H, \;\langle Q^{LH}_2(f^L, f^H), v^L \rangle = 2\pi B^{LH} n^L n^H u^L. \\
        &\langle Q^{HL}_1(f^H, f^{L}), |v^H|^2 \rangle = -4 \pi B^{HL} n^H n^L (2 T^H - 2 T^L + |u^H|^2 - |u^L|^2), 
    \end{align*}
    }
\end{lemma}
Using Lemma~\ref{lemma:collision_explicit_formula} and~\eqref{equation:conservation_properties}, we obtain a linear system for $U^{n+1}_L$ and $U^{n+1}_H$. To resolve the stiffness, we treat the terms $u^L - \varepsilon u^H$ \textit{implicitly}. Thus, we approximate the moments of $f^{n+1}_L,\, f^{n+1}_H$ as follows
{\small
\begin{equation*}
    \begin{aligned}
        \tilde{n}^{n+1}_L &= n^n_L, \quad \tilde{n}^{n+1}_H = n^n_H, \\
        \tilde{u}^{n+1}_L &= u^n_L - 2\pi B^{HL} n^n_H\Delta t \left[ \frac{1}{\tau} (\tilde{u}^{n+1}_L - \varepsilon \tilde{u}^{n+1}_H) - \frac{\varepsilon^2}{\tau} u^n_L\right], \\
        \tilde{u}^{n+1}_H &= u^n_H + 2\pi B^{HL} n^n_L \Delta t \ \frac{\varepsilon}{\tau}  (\tilde{u}^{n+1}_L - \varepsilon \tilde{u}^{n+1}_H), \\
        \tilde{T}^{n+1}_L &= T^n_L - 2\pi B^{HL} n^n_H \Delta t \left[ \frac{\varepsilon}{\tau}  (\tilde{u}^{n+1}_L - \varepsilon \tilde{u}^{n+1}_H) \cdot u^n_H 
        - \frac{\varepsilon^2}{\tau}  (2T^n_H - 2T^n_L - |u^n_L|^2) \right], \\
        \tilde{T}^{n+1}_H &= T^n_H + 2\pi B^{HL} n^n_L \Delta t \left[ \frac{\varepsilon}{\tau}  (\tilde{u}^{n+1}_L - \varepsilon \tilde{u}^{n+1}_H) \cdot u^n_H  
        - \frac{\varepsilon^2}{\tau}  (2T^n_H - 2T^n_L - |u^n_L|^2) \right].
    \end{aligned}
\end{equation*}}To distinguish the solutions of the above system from the true moments of $f^{n+1}_L,\, f^{n+1}_H$, we denote them as $\widetilde{U}^{n+1}_L = (\tilde{n}^{n+1}_L, \tilde{u}^{n+1}_L, \tilde{T}^{n+1}_L)$, $\widetilde{U}^{n+1}_H = (\tilde{n}^{n+1}_H, \tilde{u}^{n+1}_H, \tilde{T}^{n+1}_H)$. The $(n+1)$-level Maxwellians in \eqref{equation:scheme_initial} are then defined according to~\eqref{equation:maxwellian_def}. 

\subsubsection*{The penalty parameter $\nu$} In the BGK-penalization approach, $\nu_{L}$ and $\nu_{H}$ are positive constants chosen for stability, which was discussed in \cite{Filbet-Jin, Yan-Jin, Gamba-Jin-Liu}. In the case of the Boltzmann operator $Q$, one typical choice of the parameter $\nu$ is given by $\nu > Q^-$,
in which the decomposition of the Boltzmann collision operator
$Q = Q^+ - fQ^-$ is made. 
 This approach guarantees the positivity of the numerical solutions \cite{Yan-Jin}. Here, we adopt the same approach by choosing
$\nu_{L} \geq Q^{LL, -} + Q^{LH, -}_\varepsilon, \quad \nu_{H} \geq Q^{HH, -} + Q^{HL, -}_\varepsilon$. In the case of Maxwell molecules, we can choose $\nu_{L} = Q^{LL, -} + Q^{LH, -}_0 =  2 \pi B^{LH} (n_L + n_H)$ and $\nu_{H} = Q^{HH, -} + Q^{HL, -}_0 = 2 \pi B^{HL} (n_L + n_H)$.

\subsubsection*{The final AP scheme}
{
We now summarize the steps to obtain $f^{n+1}_L,\, f^{n+1}_H$ from $f^n_L,\, f^n_H$ in \textbf{Algorithm~\ref{algorithm:AP}}. First, we compute the moments $U^n_L,\, U^n_H$ of $f^n_L,\, f^n_H$ by~\eqref{equation:moments}, from which we define $\nu^{n}_L$, $\nu^{n}_H$ and $M^n_L$, $M^n_H$. Then we update the moments $\widetilde{U}^{n+1}_L$ and $\widetilde{U}^{n+1}_H$ \textit{a priori} using~\eqref{equation:moment_update}, which in turn define $\tilde{\nu}^{n+1}_L$, $\tilde{\nu}^{n+1}_H$ and $\widetilde{M}^{n+1}_L$, $\widetilde{M}^{n+1}_H$. Note that $\widetilde{M}^{n+1}_L$ and $\widetilde{M}^{n+1}_H$ are only approximations of the true Maxwellians $M^{n+1}_L$ and $M^{n+1}_H$ associated with $f^{n+1}_L$ and $f^{n+1}_H$. Finally, we use $\widetilde{M}^{n+1}_L$, $\widetilde{M}^{n+1}_H$ to replace $M^{n+1}_L$, $M^{n+1}_H$ in~\eqref{equation:scheme_initial} and compute the distribution functions at time $t^{n+1}$.}

\begin{algorithm}
    \small
    \caption{\small {The AP time discretization}}
    \label{algorithm:AP}
    \begin{algorithmic}[1]
    
    \State {\textbf{Input:} $f^n_L,\, f^n_H$}
    
    \State {Compute moments at time $t^n$:} 

    \begin{minipage}{0.8\textwidth}
        \resizebox{1\textwidth}{!}{%
          {$\displaystyle
          \begin{aligned}
            n^n_L &= \sum_{i=1}^{N_v^2} f^n_L(v_i), \;
            n^n_H = \sum_{i=1}^{N_v^2} f^n_H(v_i), \;
            u^n_L = \frac{1}{n^n_L} \sum_{i=1}^{N_v^2} f^n_L(v_i) v_i, \;
            u^n_H = \frac{1}{n^n_H} \sum_{i=1}^{N_v^2} f^n_H(v_i) v_i, \\
            T^n_L &= \frac{1}{2 n^n_L} \sum_{i=1}^{N_v^2} f^n_L(v_i) |v_i - u^n_L|^2, \;
            T^n_H = \frac{1}{2 n^n_H} \sum_{i=1}^{N_v^2} f^n_H(v_i) |v_i - u^n_H|^2.
          \end{aligned}$}
        }
    \end{minipage}

    \State {Define $\nu^{n}_L = 2 \pi B^{LH} (n^n_L + n^n_H)$ and $\nu^{n}_H = 2 \pi B^{HL} (n^n_L + n^n_H)$ and Maxwellians at time $t^n$:} 

    \begin{minipage}{0.8\textwidth}
        \resizebox{0.9\textwidth}{!}{\footnotesize
          {$\displaystyle
          \begin{aligned}
            M^n_L(v_i) = \frac{n^n_L}{2\pi T^n_L} \exp\left(-\frac{|v_i-u^n_L|^2}{2T^n_L}\right), \quad
            M^n_H(v_i) = \frac{n^n_H}{2\pi T^n_H} \exp\left(-\frac{|v_i-u^n_H|^2}{2T^n_H}\right).
          \end{aligned}
          $}
        }
      \end{minipage}

    \State {Update the moments at time $t^{n+1}$:}

    \resizebox{0.9\textwidth}{!}{%
    \begin{minipage}{1\textwidth}
        {\footnotesize
    {\begin{equation}
        \label{equation:moment_update}
    \begin{aligned}
        \tilde{n}^{n+1}_L &= n^n_L, \quad \tilde{n}^{n+1}_H = n^n_H, \\
        \tilde{u}^{n+1}_L &= u^n_L - 2\pi B^{HL} n^n_H\Delta t \left[ \frac{1}{\tau} (\tilde{u}^{n+1}_L - \varepsilon \tilde{u}^{n+1}_H) - \frac{\varepsilon^2}{\tau} u^n_L\right], \\
        \tilde{u}^{n+1}_H &= u^n_H + 2\pi B^{HL} n^n_L \Delta t \ \frac{\varepsilon}{\tau}  (\tilde{u}^{n+1}_L - \varepsilon \tilde{u}^{n+1}_H), \\
        \tilde{T}^{n+1}_L &= T^n_L - 2\pi B^{HL} n^n_H \Delta t \left[ \frac{\varepsilon}{\tau}  (\tilde{u}^{n+1}_L - \varepsilon \tilde{u}^{n+1}_H) \cdot u^n_H 
        - \frac{\varepsilon^2}{\tau}  (2T^n_H - 2T^n_L - |u^n_L|^2) \right], \\
        \tilde{T}^{n+1}_H &= T^n_H + 2\pi B^{HL} n^n_L \Delta t \left[ \frac{\varepsilon}{\tau}  (\tilde{u}^{n+1}_L - \varepsilon \tilde{u}^{n+1}_H) \cdot u^n_H  
        - \frac{\varepsilon^2}{\tau}  (2T^n_H - 2T^n_L - |u^n_L|^2) \right].
    \end{aligned}
    \end{equation}}
    }
    \end{minipage}%
    }
    \State {Define $\tilde{\nu}^{n+1}_L = 2 \pi B^{LH} (\tilde{n}^{n+1}_L + \tilde{n}^{n+1}_H)$ and $\tilde{\nu}^{n+1}_H = 2 \pi B^{HL} (\tilde{n}^{n+1}_L + \tilde{n}^{n+1}_H)$ and Maxwellians at time $t^{n+1}$:} 

    \begin{minipage}{0.8\textwidth}
        \resizebox{1\textwidth}{!}{ \footnotesize
          {$\displaystyle
          \begin{aligned}
            \widetilde M_L^{n+1}(v_i)
              = \frac{\tilde n_L^{n+1}}{2\pi \tilde T_L^{n+1}}
                 \exp\!\Bigl(-\tfrac{|v_i-\tilde u_L^{n+1}|^2}{2\tilde T_L^{n+1}}\Bigr), \quad
            \widetilde M_H^{n+1}(v_i)
              = \frac{\tilde n_H^{n+1}}{2\pi \tilde T_H^{n+1}}
                 \exp\!\Bigl(-\tfrac{|v_i-\tilde u_H^{n+1}|^2}{2\tilde T_H^{n+1}}\Bigr).
          \end{aligned}
          $}
        }
      \end{minipage}

    \State {Update distribution functions at time $t^{n+1}$:} 

    \resizebox{0.8\textwidth}{!}{%
    \begin{minipage}{1\textwidth}
        {\small
    {\begin{equation}
    \label{equation:scheme}
    \begin{aligned}
        \frac{f^{n+1}_{L}-f^{n}_{L}}{\Delta t}
        =& \frac{1}{\tau}\left(Q^{LL}(f^n_L,f^n_L) + Q^{LH}_{\text{AE}}(f^n_L,f^n_H)
        - \nu^{n}_{L}(M^n_{L}-f^{n}_{L})\right) \\ & + \frac{1}{\tau}\tilde{\nu}^{n+1}_{L}(\widetilde{M}^{n+1}_{L}-f^{n+1}_{L}), \\
        \frac{f^{n+1}_{H}-f^{n}_{H}}{\Delta t}
        =& \frac{\varepsilon}{\tau}\left(Q^{HH}(f^{n}_{H},f^{n}_{H}) + Q^{HL}_{\text{AE}}(f^n_H,f^n_L)
        -\nu^{n}_{H}(M^n_{H}-f^{n}_{H})\right) \\ & + \frac{\varepsilon}{\tau} \tilde{\nu}^{n+1}_{H} (\widetilde{M}^{n+1}_{H}-f^{n+1}_{H}).
    \end{aligned}
    \end{equation}}
    }
    \end{minipage}%
    }
    \State {\textbf{Output:} $f^{n+1}_L,\, f^{n+1}_H$}
    \end{algorithmic}
\end{algorithm}

\vspace{2mm}

\subsubsection*{Asymptotic behavior of the moment system}
We now analyze the asymptotic behavior of the scheme's moments across different time scales by examining the solutions to~\eqref{equation:moment_update}. 
For simplicity, we define:
\begin{equation*}
    \alpha = 2\pi B^{HL} n^n_H \frac{\Delta t}{\tau}, \quad 
    \beta = 2\pi B^{HL} n^n_L \frac{\varepsilon \Delta t}{\tau}, \quad 
    \gamma = 2\pi B^{HL} n^n_H \frac{\varepsilon^2 \Delta t}{\tau}.
\end{equation*}
From~\eqref{equation:moment_update}, we have the following explicit updates
{\small
\begin{equation}
    \label{equation:velocity_update}
    \begin{aligned}
        \tilde{u}^{n+1}_H &= \frac{\beta(u^n_L + \gamma u^n_L) + (1 + \alpha)u^n_H}{1 + \alpha + \beta\varepsilon}, \quad
        \tilde{u}^{n+1}_L - \varepsilon \tilde{u}^{n+1}_H = \frac{u^n_L + \gamma u^n_L - \varepsilon u^n_H}{1 + \alpha + \beta\varepsilon}
    \end{aligned}
\end{equation}
}Based on this, we can summarize the asymptotic behavior of the scheme's moments across different time scales in Table~\ref{tab:epochal}.
\begin{table}[ht]
    \centering
    \resizebox{1\textwidth}{!}{%
    \begin{tabular}{|p{2.48cm}|c|c|c|}
    \hline 
    $\begin{array}{l}
         \; \\[-5pt]
    \textbf{\small Time Scale} \\ [-5pt]
    \;
    \end{array}$
    & \textbf{Fastest} ($\tau = 1$) 
    & \textbf{Intermediate} ($\tau = \varepsilon$) 
    & \textbf{Slowest} ($\tau = \varepsilon^2$) \\
    \hline
    %% Parameters
    $\begin{array}{l}
        \; \\[-5pt]
    \textbf{\small Parameters} \\ [-5pt]
    \;
    \end{array}$
    & 
    $\begin{array}{l}
    \alpha = \mathcal{O}(\Delta t), \
    \beta = \mathcal{O}(\Delta t \varepsilon), \
    \gamma = \mathcal{O}(\Delta t \varepsilon^2) \\
    \end{array}$
    & 
    $\begin{array}{l}
    \alpha = \mathcal{O}(\Delta t \varepsilon^{-1}), \
    \beta = \mathcal{O}(\Delta t), \
    \gamma = \mathcal{O}(\Delta t \varepsilon)
    \end{array}$ & 
    $\begin{array}{l}
    \alpha = \mathcal{O}(\Delta t \varepsilon^{-2}), \
    \beta = \mathcal{O}(\Delta t \varepsilon^{-1}), \
    \gamma = \mathcal{O}(\Delta t)
    \end{array}$ \\
    \hline
    %% Velocity
    $\begin{array}{l}
        \; \\
    \textbf{\small Velocity} \\ [-5pt]
    \;
    \end{array}$
    & 
    $\begin{array}{l}
    \tilde{u}^{n+1}_H = u^n_H + \mathcal{O}(\Delta t \varepsilon) \\[3pt]
    \tilde{u}^{n+1}_L - \varepsilon \tilde{u}^{n+1}_H = \frac{u^n_L - \varepsilon u^n_H}{1 + \alpha} + \mathcal{O}(\Delta t \varepsilon^2) 
    \end{array}$ & 
    $\begin{array}{l}
    \tilde{u}^{n+1}_H = u^n_H + \varepsilon u^n_L + \mathcal{O}(\Delta t \varepsilon^2) \\[3pt]
    \tilde{u}^{n+1}_L - \varepsilon \tilde{u}^{n+1}_H = \mathcal{O}(\Delta t \varepsilon) u^n_L + \mathcal{O}(\Delta t \varepsilon^2)
    \end{array}$ & 
    $\begin{array}{l}
    \tilde{u}^{n+1}_H = u^n_H + \mathcal{O}(\Delta t \varepsilon^2) \\[3pt]
    \tilde{u}^{n+1}_L - \varepsilon \tilde{u}^{n+1}_H = \mathcal{O}(\Delta t \varepsilon^2)
    \end{array}$ \\
    \hline
    %% Temperature
    $\begin{array}{l}
        \; \\
    \textbf{\small Temperature} \\ 
    \;
    \end{array}$
    & 
    $\begin{array}{l}
    \tilde{T}^{n+1}_L = T^n_L + \mathcal{O}(\Delta t \varepsilon) \\[3pt]
    \tilde{T}^{n+1}_H = T^n_H + \mathcal{O}(\Delta t \varepsilon)
    \end{array}$ & 
    $\begin{array}{l}
    \tilde{T}^{n+1}_L = T^n_L + \mathcal{O}(\Delta t \varepsilon) \\[3pt]
    \tilde{T}^{n+1}_H = T^n_H + \mathcal{O}(\Delta t \varepsilon)
    \end{array}$ & 
    $\begin{array}{l}
    \tilde{T}^{n+1}_L = T^n_L + 2\pi B^{HL} n^n_H \Delta t (2T^n_H - 2T^n_L) + \mathcal{O}(\Delta t \varepsilon) \\[3pt]
    \tilde{T}^{n+1}_H = T^n_H - 2\pi B^{HL} n^n_L \Delta t (2T^n_H - 2T^n_L) + \mathcal{O}(\Delta t \varepsilon)
    \end{array}$ \\
    \hline
    %% Behavior
    \textbf{\small $\mathcal{O}(1)$ behavior} & 
    \begin{minipage}{0.45\textwidth}
    \vspace{2mm}
    $\tilde{u}_L$ relaxes to $\varepsilon \tilde{u}_H$ with rate $\mathcal{O}(\frac{1}{1 + \alpha})$; \newline
    $\tilde{T}_L$ and $\tilde{T}_H$ remain constant.
    \vspace{2mm}
    \end{minipage}
    % Velocity relaxation
    & 
    \begin{minipage}{0.5\textwidth}
    $\tilde{u}_L = \varepsilon \tilde{u}_H$; 
    $\tilde{T}_L$ and $\tilde{T}_H$ remain constant.
    \end{minipage} & 
    \begin{minipage}{0.5\textwidth}
    % Temperature relaxation
    $\tilde{u}_L = \varepsilon \tilde{u}_H$; 
    $\tilde{T}_L$ and $\tilde{T}_H$ relax towards each other according to~\eqref{equation:macro}.
    \end{minipage} \\[1pt]
    \hline
    \end{tabular}
    }
    \caption{\footnotesize Epochal relaxation of velocity and temperature across three time scales.}
    \label{tab:epochal}
\end{table}

The moment updating process exhibits velocity relaxation of $u^L$ towards $\varepsilon u^H$ at the fastest time scale. One recalls from Remark~\ref{remark:heavy_species_macros} that the mean velocity of the heavy species is $\varepsilon u^H$. 

From ~\eqref{equation:moment_update} we also see that the temperature relaxation similar to~\eqref{equation:macro} happens at the slowest time scale \textit{only} after the velocity relation $u^L - \varepsilon u^H = \mathcal{O}(\varepsilon^2)$ is reached. This requires us to treat the term $u^L - \varepsilon u^H$ implicitly, instead of just $u^L$.

\subsubsection*{Consistency of the moment update}
We now show that the implicit treatment of the velocities in the moment update ~\eqref{equation:moment_update} is consistent with the scheme~\eqref{equation:scheme}
\begin{proposition}
    The moment update~\eqref{equation:moment_update} is consistent with the scheme~\eqref{equation:scheme} meaning that we have the following error bounds
    {\small
    \begin{equation*}
        \widetilde{U}^{n+1}_L - U^{n+1}_L = \mathcal{O}\left(\frac{\Delta t^3 }{\tau + \nu^{n+1}_L \Delta t} \right), \quad \widetilde{U}^{n+1}_H - U^{n+1}_H = \mathcal{O}\left(\frac{\varepsilon\Delta t^3 }{\tau + \nu^{n+1}_H \Delta t} \right),
    \end{equation*}
    }where $\widetilde{U}^{n+1}_L$ and $\widetilde{U}^{n+1}_H$ are the solution of the system~\eqref{equation:moment_update}, while $U^{n+1}_L$ and $U^{n+1}_H$ are the moments of $f^{n+1}_L$ and $f^{n+1}_H$ computed from directly taking moments of the scheme~\eqref{equation:scheme}.
    \begin{proof}
        \small
        We only prove for the light species equation, and the proof for the heavy species is similar. By taking moments of the scheme~\eqref{equation:scheme}, we have
        {\small
        \begin{equation}
        \label{eqn:scheme1}
            \begin{aligned}
                \frac{U^{n+1}_L - U^n_L}{\Delta t}
                &= \frac{1}{\tau} \langle Q^{LH}_0(f^n_L,f^n_H), \phi(v^L) \rangle + \frac{\varepsilon}{\tau} \langle Q^{LH}_1(f^n_L ,f^n_H), \phi(v^L) \rangle \\[5pt]
                & + \frac{\varepsilon^2}{\tau} \langle Q^{LH}_2(f^n_L,f^n_H), \phi(v^L) \rangle  + \frac{1}{\tau} \nu^{n+1}_L (\widetilde{U}^{n+1}_L - U^{n+1}_L), 
            \end{aligned}
        \end{equation}
        }where $\phi(v) = 1,\ v, \ |v|^2$.
        Comparing with equations ~\eqref{equation:moment_update}, one can rewrite the equation \eqref{eqn:scheme1} as
        {\small
        \begin{align*}
            \frac{U^{n+1}_L - U^n_L}{\Delta t}
            &= \frac{\widetilde{U}^{n+1}_L - U^n_L}{\Delta t} +  \frac{1}{\tau} \nu^{n+1}_L (\widetilde{U}^{n+1}_L - U^{n+1}_L) \\ 
            &+ \mathcal{O}(\frac{\Delta t}{\tau}) (\tilde{u}^{n+1}_L - u^n_L)
             + \mathcal{O}(\frac{\varepsilon \Delta t}{\tau})(\tilde{u}^{n+1}_H - u^n_H).
        \end{align*}
        }This gives
        {\small
        \begin{equation*}
            \begin{aligned}
                U^{n+1}_L - \widetilde{U}^{n+1}_L
                &= \frac{\Delta t }{\tau + \nu^{n+1}_L \Delta t} \mathcal{O}(\Delta t) \left[
                 (\tilde{u}^{n+1}_L - u^n_L)
                + \varepsilon (\tilde{u}^{n+1}_H - u^n_H)
                \right].
            \end{aligned}
        \end{equation*}
        }We are reduced to bounding the magnitude of the one-step changes $\tilde{u}^{n+1}_L - u^n_L$ and $\tilde{u}^{n+1}_H - u^n_H$. From~\eqref{equation:velocity_update}, we compute
        {\small
        \begin{align*}
            \tilde{u}^{n+1}_L - u^n_L = -\alpha W^{n+1} + \gamma u^n_L, \quad
            \tilde{u}^{n+1}_H - u^n_H = \beta W^{n+1}, 
        \end{align*}
        }where $W^{n+1} = \tilde{u}^{n+1}_L - \varepsilon \tilde{u}^{n+1}_H$.
        Using the conclusion from Table~\ref{tab:epochal}, we have the following case studies
        {\small
        \begin{itemize}
            \item $\tau = 1$:
            $\alpha = \mathcal{O}(\Delta t)$, $\beta = \mathcal{O}(\varepsilon\Delta t)$, $\gamma = \mathcal{O}(\varepsilon^2\Delta t)$, $W^{n+1} = \mathcal{O}(1)$
            \vspace{5pt}
            \item $\tau = \varepsilon$:
            $\alpha = \mathcal{O}(\frac{\Delta t}{\varepsilon})$, $\beta = \mathcal{O}(\Delta t)$, $\gamma = \mathcal{O}(\varepsilon\Delta t)$, $W^{n+1} = \mathcal{O}(\varepsilon)$
            \vspace{5pt}
            \item $\tau = \varepsilon^2$:
            $\alpha = \mathcal{O}(\frac{\Delta t}{\varepsilon^2})$, $\beta = \mathcal{O}(\frac{\Delta t}{\varepsilon})$, $\gamma = \mathcal{O}(\Delta t)$, $W^{n+1} = \mathcal{O}(\varepsilon^2)$
        \end{itemize}
        \vspace{5pt}
        }Therefore, $\tilde{u}^{n+1}_L - u^n_L = \mathcal{O}(\Delta t), \quad \tilde{u}^{n+1}_H - u^n_H = \mathcal{O}(\Delta t)$. 
        This completes the proof.
    \end{proof}
\end{proposition}

% \begin{remark}
%     For higher-order IMEX-RK methods of the scheme~\eqref{equation:scheme}, the moment update~\eqref{equation:moment_update} should also be made IMEX-RK to preserve the order.
% \end{remark}

\begin{remark}
    \label{remark:collision_kernel}
    For general collision kernels, explicit formulas, such as those provided in Lemma~\ref{lemma:collision_explicit_formula}, are typically unavailable or difficult to derive. Hence, implicit treatment of the stiff terms in~\eqref{equation:taking_moments} is not easy. Exploring more efficient moment update methods for general collision kernels is deferred to future research.
\end{remark}

\subsubsection*{The AP property}
In this part, we prove the AP property of the scheme~\eqref{equation:moment_update}-\eqref{equation:scheme}. In particular, as $\tau$ varies across three time scales, our scheme can capture the epochal relaxation process in subsection~\ref{subsection:epochal_relaxation} at the kinetic level.
\begin{theorem}
    \label{theorem:AP}
   The numerical solutions $f^n_L$, $f^n_H$ given by~\eqref{equation:moment_update}-\eqref{equation:scheme} satisfy the following asymptotic properties at different time scales:
    \begin{itemize}
        \item At the fastest time scale, $f^n_L$ relaxes toward a centered Maxwellian; 
        \item At the intermediate time scale, if $f^n_L - M^n_L = \mathcal{O}(\varepsilon)$, then $f^{n+1}_L - M^{n+1}_L = \mathcal{O}(\varepsilon)$. Moreover, the numerical solution $f^n_H$ relaxes toward a Maxwellian; 
        \item At the slowest time scale, if $f^n_L - M^n_L = \mathcal{O}(\varepsilon)$ and $f^n_H - M^n_H = \mathcal{O}(\varepsilon)$, then $f^{n+1}_L - M^{n+1}_L = \mathcal{O}(\varepsilon)$ and $f^{n+1}_H - M^{n+1}_H = \mathcal{O}(\varepsilon)$. In particular, as $\varepsilon \rightarrow 0$, 
        the numerical scheme \eqref{equation:moment_update}-\eqref{equation:scheme} automatically becomes a consistent discretization of the macroscopic limit equation \eqref{equation:macro}.
    \end{itemize}
    \begin{proof}
        \small
        \hfill
        \begin{itemize}
            \item Fastest time scale ($\tau = 1$). 
            {Setting $\tau = 1$ in~\eqref{equation:scheme}, one has
            {\small
            \begin{equation}
            \label{equation:scheme_tau1}
                \begin{aligned}
                    \frac{f^{n+1}_{L}-f^{n}_{L}}{\Delta t}
                    =& \left(Q^{LL}(f^n_L,f^n_L) + Q^{LH}_{\text{AE}}(f^n_L,f^n_H)
                    - \nu^{n}_{L}(M^n_{L}-f^{n}_{L})\right) \\
                    & + \tilde{\nu}^{n+1}_{L}(\widetilde{M}^{n+1}_{L}-f^{n+1}_{L}), \\[4pt]
                    \frac{f^{n+1}_{H}-f^{n}_{H}}{\Delta t}
                    =& \varepsilon\left(Q^{HH}(f^{n}_{H},f^{n}_{H}) + Q^{HL}_{\text{AE}}(f^n_H,f^n_L)
                    -\nu^{n}_{H}(M^n_{H}-f^{n}_{H})\right) \\
                    & + \varepsilon \tilde{\nu}^{n+1}_{H} (M^{n+1}_{H}-f^{n+1}_{H}).
                \end{aligned}
            \end{equation}
            }From~\eqref{equation:AE}, we have $Q^{LH}_{\text{AE}} = Q^{LH}_0 + \mathcal{O}(\varepsilon)$. Hence, the scheme~\eqref{equation:scheme_tau1} can be rewritten as
            {\small
            \begin{align*}
                \frac{f^{n+1}_{L}-f^{n}_{L}}{\Delta t}
                =& Q^{LL}(f^n_L,f^n_L) + Q^{LH}_0(f^n_L,f^n_H) + \mathcal{O}(\varepsilon) \\
                &- \nu^{n}_{L}(M^n_{L}-f^{n}_{L}) + \tilde{\nu}^{n+1}_{L}(\widetilde{M}^{n+1}_{L}-f^{n+1}_{L}) \\
                =& Q^{LL}(f^n_L,f^n_L) + Q^{LH}_0(f^n_L,f^n_H) + \mathcal{O}(\Delta t) + \mathcal{O}(\varepsilon), \\[4pt]
                \frac{f^{n+1}_{H}-f^{n}_{H}}{\Delta t}
                =& \mathcal{O}(\varepsilon).
            \end{align*}}}
            It is clear that the scheme is a consistent discretization of the equations
            {\small
            \begin{equation}
            \label{equation:fL_leadingOrder}
                \partial_t f^L = Q^{LL}(f^L, f^L) + Q^{LH}_0(f^L, f^H), \quad \partial_t f^H = 0.
            \end{equation}
            }up to an error of order $\mathcal{O}(\varepsilon)$. 
            {
            Hence, based on \cite[Proposition 5.10, Corollary 5.11]{Degond-Lucquin-Desreux}, the numerical solution $f_L^n$ relaxes towards a centered Maxwellian. }
            \item Intermediate time scale ($\tau = \varepsilon$). We rewrite the first equation in~\eqref{equation:scheme} as
            \begin{equation*}
                \begin{aligned}
                    f^{n+1}_L = 
                    & \frac{\nu^{n+1}_{L} \Delta t}{\varepsilon + \nu^{n+1}_{L} \Delta t} M^{n+1}_L 
                    + \frac{\varepsilon}{\varepsilon + \nu^{n+1}_{L} \Delta t} f^n_L \\
                    & + \frac{\Delta t}{\varepsilon + \nu^{n+1}_{L} \Delta t} 
                    \left( Q^{LL}(f^n_L,f^n_L) + Q^{LH}_{\text{AE}}(f^n_L,f^n_H) - \nu^{n}_{L} (M^l_L - f^n_L) \right),
                \end{aligned}
            \end{equation*}
            {The collision operator $Q^{LL}$ is bilinear and vanishes for Maxwellian distributions, i.e., $Q^{LL}(M,M)=0$ (see (6.2) and (7.2) in Chapter 2 in \cite{Cercignani}). 
            By assumption, $f^n_L = M^n_L + \mathcal{O}(\varepsilon)$. 
            Actually, if we denote $g_L^n = \frac{f_L^n - \mathcal{M}_L^n}{\varepsilon}$, then we have $g_L^n = \mathcal{O}(1)$ and 
            it follows from the bilinearity of $Q^{LL}$ that
            \begin{align*}
            Q^{LL}(f^n_L, f^n_L) &= Q^{LL}(M^n_L + \varepsilon g^n_L, M^n_L + \varepsilon g^n_L) \\
            &= Q^{LL}(M^n_L, M^n_L) + \mathcal{O}(\varepsilon) = \mathcal{O}(\varepsilon).
            \end{align*}}Similarly, combining the following three facts:
            \vspace*{2pt}
            \begin{itemize}
            \item $f^n_L - M^n_L = \mathcal{O}(\varepsilon)$;
            \vspace*{2pt}
            \item $M^n_L - M^n_{L, 0} = \mathcal{O}(\varepsilon)$ after the first time step (see Table~\ref{tab:epochal}), where $M^n_{L,0}$ is the Maxwellian that shares the density, temperature with $M^n_L$ but with zero mean velocity; 
            \vspace*{2pt}
            \item $Q^{LH}_0(M^L_0, f^n_H) = 0$ by  Property 1 in Proposition~\ref{prop:asymptotic_expansions},
            \end{itemize}
            we can show $Q^{LH}_0(f^n_L, f^n_H) = \mathcal{O}(\varepsilon)$.
            Finally, 
            $f^{n+1}_L - M^{n+1}_L = \mathcal{O}(\varepsilon)$.
            By a similar argument as in the fast time scale and $Q^{HL}_0(f^n_H, M^n_{L, 0}) = 0$, we can show that the scheme for $f_H$ is consistent with
            \begin{equation}
            \label{equation:fH_leadingOrder}
                \partial_t f^H = Q^{HH}(f^H, f^H), 
            \end{equation}
            with an error of $\mathcal{O}(\varepsilon)$. 
            This implies that the collision process of $f^H$ subject to $Q^{HH}$ happens at this time scale.
            \item Slowest time scale: $\tau = \varepsilon^2$.
            At this time scale, we have
            {\small
            \begin{equation*}
                \begin{aligned}
                    f^{n+1}_L = 
                    & \frac{\nu^{n+1}_{L} \Delta t}{\varepsilon^2 + \nu^{n+1}_{L} \Delta t} M^{n+1}_L 
                    + \frac{\varepsilon^2}{\varepsilon^2 + \nu^{n+1}_{L} \Delta t} f^n_L \\
                    & + \frac{\Delta t}{\varepsilon^2 + \nu^{n+1}_{L} \Delta t} 
                    \left( Q^{LL}(f^n_L,f^n_L) + Q^{LH}_{\text{AE}}(f^n_L,f^n_H) - \nu^{n}_{L} (M^l_L - f^n_L) \right), \\
                    f^{n+1}_H = 
                    & \frac{\nu^{n+1}_{H} \Delta t}{\varepsilon + \nu^{n+1}_{H} \Delta t} M^{n+1}_H
                    + \frac{\varepsilon}{\varepsilon + \nu^{n+1}_{H} \Delta t} f^n_H \\
                    & + \frac{\Delta t}{\varepsilon + \nu^{n+1}_{H} \Delta t} 
                    \left( Q^{HH}(f^n_H,f^n_H) + Q^{HL}_{\text{AE}}(f^n_H,f^n_L) - \nu^{n}_{H} (M^l_H - f^n_H \right).
                \end{aligned}
            \end{equation*}
            }
        \end{itemize}
        By the same argument as for the intermediate time scale, we can derive 
        \begin{equation*}
            f^{n+1}_L - M^{n+1}_L = \mathcal{O}(\varepsilon), \quad f^{n+1}_H - M^{n+1}_H = \mathcal{O}(\varepsilon). 
        \end{equation*} 
    \end{proof}
\end{theorem}

\section{Numerical Examples}
\label{sec:numerics}

In this section, we present several numerical examples
to demonstrate the accuracy and efficiency of our numerical schemes. 
We consider two-dimensional problem in velocity, with the computational domain $v\in [-20, 20]^2$. 
For simplicity, collision kernels are given as 
$ 
B^{LL} = B^{HH} = \frac{1}{4\pi}, \;
B^{LH} = B^{HL} = \frac{1}{8\pi}.
$
The collision frequency in~\eqref{equation:macro} becomes
$ 
    \lambda(T^L) = 2\pi B^{HL}T^L.
$
Let the initial distributions be the double-peak Maxwellians: 
{\small
\begin{equation*}
    f^L = \frac{1}{2}\mathcal{M}_{n^L, u^L_1, T^L} + \frac{1}{2}\mathcal{M}_{n^L, u^L_2, T^L}, \quad 
    f^H = \frac{1}{2}\mathcal{M}_{n^H, u^H_1, T^H} + \frac{1}{2}\mathcal{M}_{n^H, u^H_2, T^H},
\end{equation*}
}with Maxwellians defined by ~\eqref{equation:maxwellian} and
{\small
\begin{gather*}
    n^L = n^H = 1, \quad
    T^L = 3, \quad T^H = 0.5, \\
    u^L_1 = (1.2,0)^\top, \quad 
    u^L_2 = -(0.5, 0)^\top, \quad 
    u^H_1 = -(1.2,0)^\top, \quad 
    u^H_2 = (0.5, 0)^\top.
\end{gather*}
}

\subsection{Approximation of the inter-particle collision operators}

In this test, we assess the \textbf{AE} method and the \textbf{SP} method as candidates for approximating inter-particle collision operators in our problem. We solve the evolution problem~\eqref{equation:main_equation} for $\tau=1$ by forward Euler method 
using the \textbf{AE} method and the \textbf{SP} method to approximate the collision operators $Q^{LH}_{\varepsilon}$ and $Q^{HL}_{\varepsilon}$.

\vspace{1em}

\noindent\textbf{Convergence study of the \textbf{AE} and \textbf{SP} methods. }
We first study the convergence in the velocity variable of the numerical solutions computed by the \textbf{AE} and \textbf{SP} methods. 
In view of the discussion in Section~\ref{subsec:AE_applicability}, 
the two methods operate differently in $\varepsilon$ regimes: the \textbf{AE} method is expected to be accurate only for small $\varepsilon$, while the \textbf{SP} method requires an increasingly larger $N_v$ as $\varepsilon$ decreases.
Therefore, for the \textbf{AE} method, we study four cases of gas mixtures
characterized by $\varepsilon = 0.2$, $\varepsilon = 0.1$, $\varepsilon = 0.03$, and $\varepsilon = 0.01$ with velocity discretization by $N_v=30$, $60$, $120$, $240$, and $480$ for all $\varepsilon$ values.
For the \textbf{SP} method, we study four cases of gas mixtures characterized by  $\varepsilon = 0.1$, $\varepsilon = 0.2$, $\varepsilon = 0.5$, and $\varepsilon = 1$. According to Proposition \ref{Prop:SP}, for each $\varepsilon$, we need to use velocity discretizations proportional to $1/\varepsilon$. The details of $\varepsilon$ and $N_v$ for the \textbf{SP} method are summarized in Table~\ref{table:parametric_study_spec}. We set the final computational time to be $t=0.5$ for all the tests.
\begin{table}
    \centering
    \resizebox{4cm}{!}{
    \begin{tabular}{|c|c|}
        \hline
        \text{$\varepsilon$} & \textbf{$N_v$} \\
        \hline
        $0.1$  & $50$, $100$, $200$, $400$, $800$\\
        $0.2$  & $50$, $100$, $200$, $400$, $800$ \\
        $0.5$  & $20$, $40$, $80$, $160$, $320$ \\
        $1$    & $10$, $20$, $40$, $80$, $160$ \\
        \hline
    \end{tabular}}
    \caption{\footnotesize Parameter settings for the \textbf{SP} method; velocity discretization is proportional to $1/\varepsilon$.}
    \label{table:parametric_study_spec}
\end{table}

\begin{figure}
    \centering
    \begin{subfigure}{0.25\textwidth}
        \centering
        \includegraphics[width=1\textwidth]{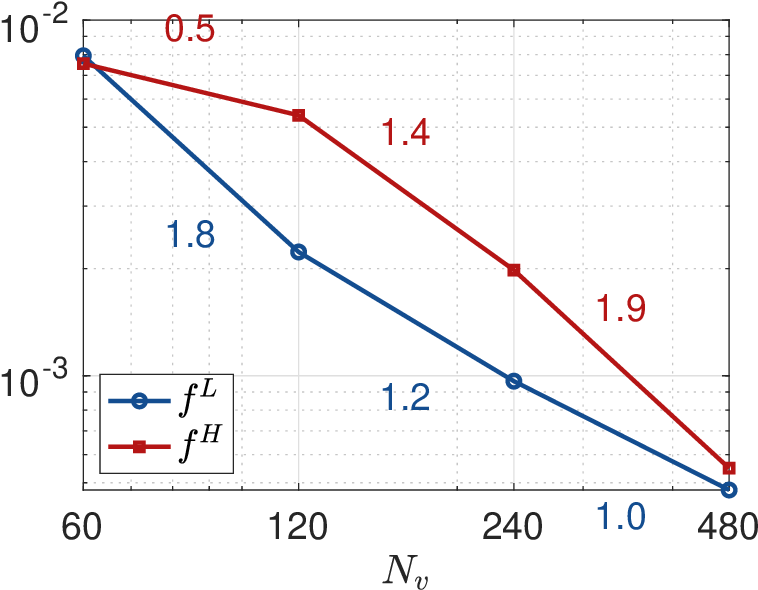}
        \caption*{\footnotesize \textbf{AE} ($\varepsilon = 0.2$)}
    \end{subfigure}
    \hspace{-6pt}
    \begin{subfigure}{0.25\textwidth}
        \centering
        \includegraphics[width=1\textwidth]{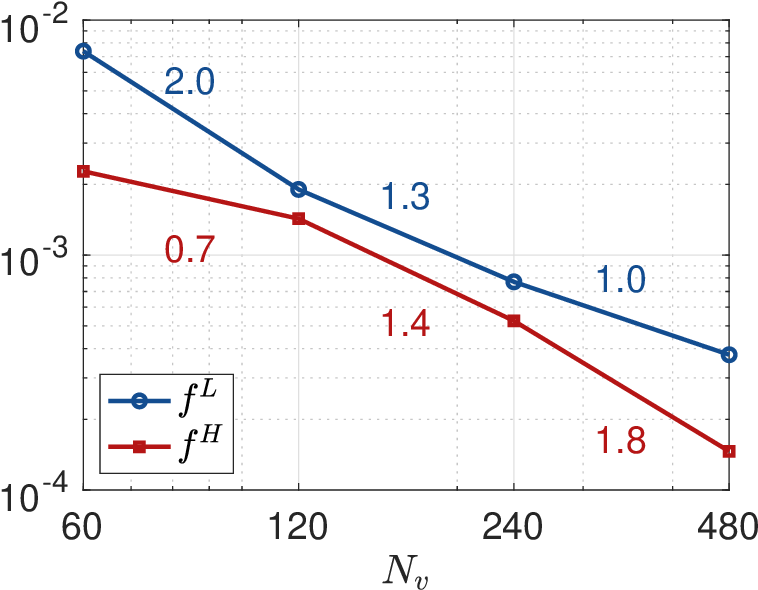}
        \caption*{\footnotesize \textbf{AE} ($\varepsilon = 0.1$)}
    \end{subfigure}
    \hspace{-6pt}
    \begin{subfigure}{0.25\textwidth}
        \centering
        \includegraphics[width=1\textwidth]{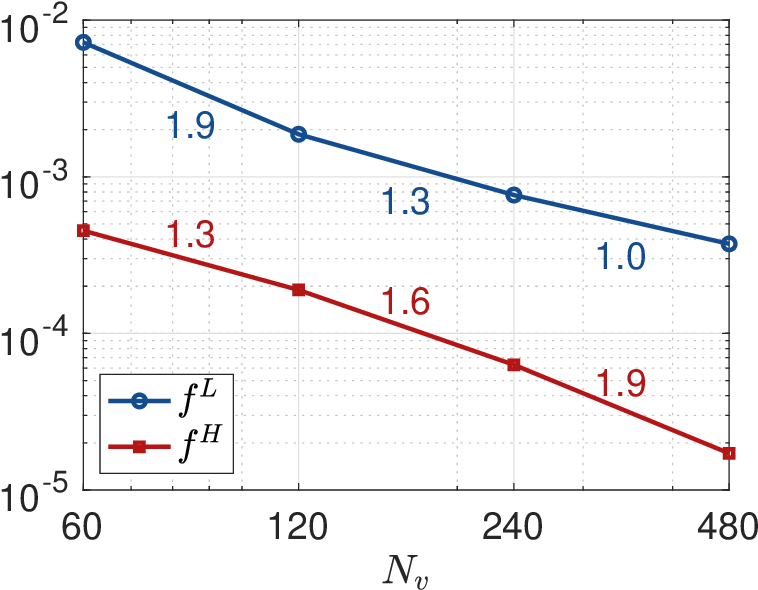}
        \caption*{\footnotesize \textbf{AE} ($\varepsilon = 0.03$)}
    \end{subfigure}
    \hspace{-6pt}
    \begin{subfigure}{0.25\textwidth}
        \centering
        \includegraphics[width=1\textwidth]{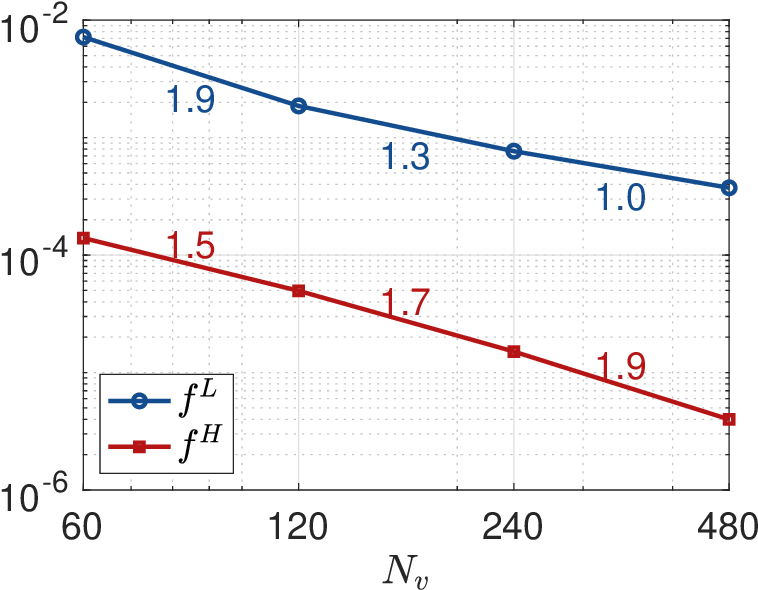}
        \caption*{\footnotesize \textbf{AE} ($\varepsilon = 0.01$)}
    \end{subfigure}
    \vskip\baselineskip
    \vspace{-1em}
    \begin{subfigure}{0.25\textwidth}
        \centering
        \includegraphics[width=1\textwidth]{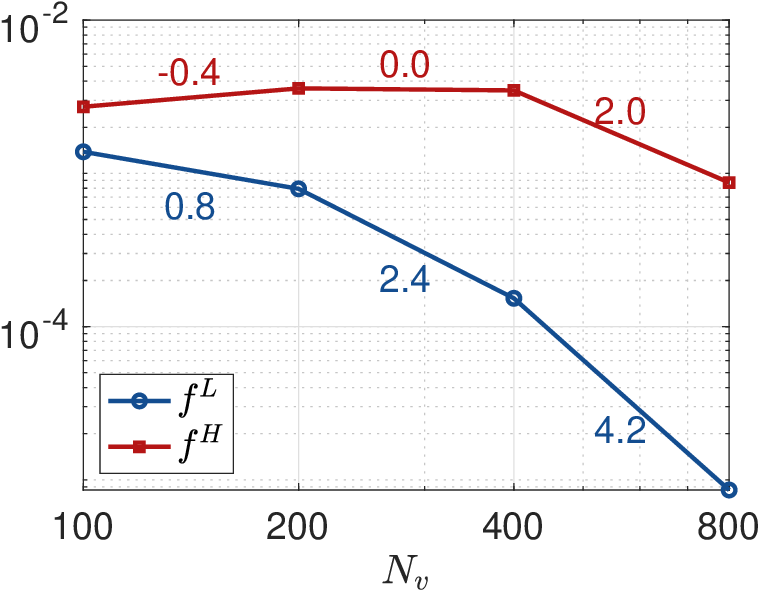}
        \caption*{\footnotesize \textbf{SP} ($\varepsilon = 0.1$)}
    \end{subfigure}
    \hspace{-6pt}
    \begin{subfigure}{0.25\textwidth}
        \centering
        \includegraphics[width=1\textwidth]{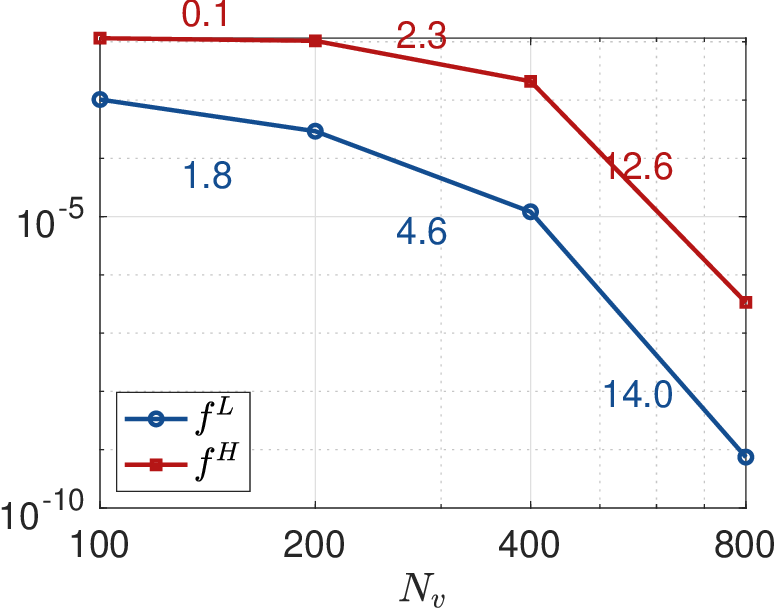}
        \caption*{\footnotesize \textbf{SP} ($\varepsilon = 0.2$)}
    \end{subfigure}
    \hspace{-6pt}
    \begin{subfigure}{0.25\textwidth}
        \centering
        \includegraphics[width=1\textwidth]{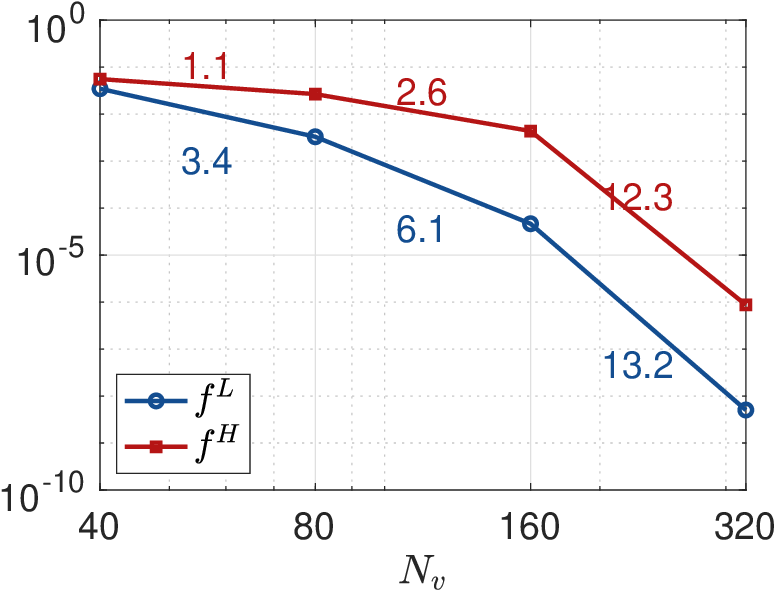}
        \caption*{\footnotesize \textbf{SP} ($\varepsilon = 0.5$)}
    \end{subfigure}
        \hspace{-6pt}
    \begin{subfigure}{0.25\textwidth}
        \centering
        \includegraphics[width=1\textwidth]{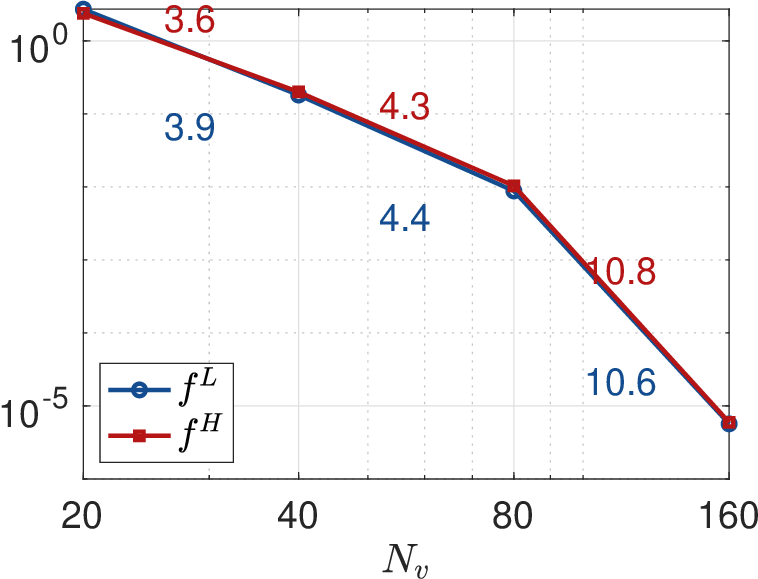}
        \caption*{\footnotesize \textbf{SP} ($\varepsilon = 1$)}
    \end{subfigure}
    \caption{\footnotesize Convergence in velocity of the \textbf{AE} and \textbf{SP} methods for different $\varepsilon$. For the \textbf{SP} method, velocity discretizations proportional to $1/\varepsilon$ are used. The y-axis is in log scale.}
    \label{figure:AE_SP_comparison}
\end{figure}

In Figure~\ref{figure:AE_SP_comparison}, we present the log-plot of the relative $\ell_v^2$ errors of numerical solutions defined by 
{\small
\begin{equation}
\label{equation:l2error}
    e_{Nv}(f) = \frac{\|f_{N_v} - f_{N_v/2}\|_{\ell_v^2}}{\|f_{N_v}\|_{\ell_v^2}}. 
\end{equation}
}Here $f_{N_v}$ denotes the numerical solution computed on a velocity grid of size $N_v \times N_v$, the slopes of each segment are plotted in the figure. We also record the computational time and orders of convergence for cases $\varepsilon = 0.1$ and $\varepsilon = 0.2$ in Table~\ref{table:time}. It is known that a numerical scheme is $k$-th order if $e_{N_v}(f) \leq \frac{C}{N_v^k}$ for $N_v \gg 1$.
In Figure~\ref{figure:AE_SP_comparison}, for the \textbf{SP} method, even though it is spectrally accurate, we observe that to maintain the same level of accuracy as $\varepsilon$ decreases, a velocity discretization proportional to $1/\varepsilon$ is necessary, which is consistent with our theoretical findings in Proposition \ref{Prop:SP}. On the other hand, from Table~\ref{table:time} we observe that, to achieve comparable accuracy for $\varepsilon = 0.1$ to that for $\varepsilon = 0.2$, the \textbf{SP} method requires a prohibitively growing computational cost. 

On the contrary, the \textbf{AE} method exhibits two distinct properties. First, there is a slight deterioration of accuracy as $\varepsilon$ increases, which coincides with our expectations, since according to the truncation in~\eqref{equation:AE}, the \textbf{AE} method fails to produce accurate results when $\varepsilon$ is large. 
Second, the method exhibits similar order of accuracy as shown in the slope pairs
$(f^L, f^H)$, by using the same $N_v$ values in different regimes, indicating that the method demonstrates a {\it uniform} accuracy especially for $\varepsilon<0.1$. From Table~\ref{table:time}, we can see that the \textbf{AE} method has consistent computational cost across different regimes. This behaviour of the numerical solutions benefits from the design of the \textbf{AE} method where $\varepsilon$ is decoupled from the collision operators. Lastly, in regimes with smaller $\varepsilon$, the \textbf{AE} method exhibits first-order accuracy for the distribution $f^L$ and second-order accuracy for the distribution $f^H$.

We summarize our findings based on the above observations. The \textbf{SP} method achieves optimal performance for mixtures with nearly equal masses ($\varepsilon \approx 1$), whereas becomes computationally infeasible in disparate mass regimes due to almost not practical requirements in velocity discretization.
On the contrary, the \textbf{AE} method demonstrates a \textit{uniformly cheap} computational cost and works consistently efficient in disparate mass regimes where $\varepsilon \ll 1$ while maintaining the accuracy. This conclusion matches with our observation in Figure~\ref{figure:AE_SP_regimes}, making the \textbf{AE} method a \textit{suitable candidate} for gas mixture simulations in disparate mass regimes (when $\varepsilon \ll 1$).
\begin{table}[hbt!]
    \centering
    \begin{minipage}[t]{0.45\textwidth}
        \centering
        \resizebox{\textwidth}{!}{
        \begin{tabular}{c | c c | c c}
        \toprule
         & \multicolumn{4}{c}{\textbf{SP} method} \\
        \cmidrule(lr){2-5}
         & \multicolumn{2}{c|}{\small Time (s)} & \multicolumn{2}{c}{\small Slope pairs ($f^L$, $f^H$)} \\
        \cmidrule(lr){2-3} \cmidrule(lr){4-5}
        $N_v$ & $\varepsilon = 0.2$ & $\varepsilon = 0.1$ & $\varepsilon = 0.2$ & $\varepsilon = 0.1$ \\
        \midrule
        200 & 1351.10 & 1298.88 & (0.1, 1.8) & (-0.4, 0.8) \\
        400 & 5695.96 & 5326.12 & (2.3, 4.6) & (0.0, 2.4) \\
        800 & 33854.44 & 33221.87 & (12.6, 14.0) & (2.0, 4.2) \\
        \bottomrule
        \end{tabular}
        }
    \end{minipage}
    \hfill
    \begin{minipage}[t]{0.4\textwidth}
        \centering
        \resizebox{\textwidth}{!}{
        \begin{tabular}{c | c c | c c}
        \toprule
         & \multicolumn{4}{c}{\textbf{AE} method} \\
        \cmidrule(lr){2-5}
         & \multicolumn{2}{c|}{\small Time (s)} & \multicolumn{2}{c}{\small Slope pairs ($f^L$, $f^H$)} \\
        \cmidrule(lr){2-3} \cmidrule(lr){4-5}
        $N_v$ & $\varepsilon = 0.2$ & $\varepsilon = 0.1$ & $\varepsilon = 0.2$ & $\varepsilon = 0.1$ \\
        \midrule
        120 & 4.32 & 4.08 & (1.8, 0.8) & (2.0, 0.7) \\
        240 & 12.31 & 12.34 &(1.2, 1.4) & (1.3, 1.4) \\
        480 & 47.37 & 46.39 & (1.0, 1.9) & (1.0, 1.8) \\
        \bottomrule
        \end{tabular}
        }
    \end{minipage}
    \caption{\footnotesize Computational times (in seconds) and accuracy (slope pairs) of the \textbf{SP} method (left) and the \textbf{AE} method (right) with different numbers of velocity points.}
    \label{table:time}
\end{table}

\noindent\textbf{Accuracy of the truncation approach in \textbf{AE}. }
The \textbf{AE} method uses truncation in $\varepsilon$ to approximate the inter-particle collision operators.
In this test, we investigate if the truncation approach is accurate, with the \textbf{SP} method as the reference solution.
We examine three \textit{moderate} mass disparity regimes $\varepsilon = 0.2$, $\varepsilon = 0.1$, and $\varepsilon = 0.05$ corresponding to mass ratios $m_H/m_L = 25$, $100$, and $400$, respectively. 
For each value of $\varepsilon$, we conduct simulations using the \textbf{AE} method with increasing velocity grid resolutions ($N_v = 20, 40, 80, 160, 320$). 
The reference solution computed by the \textbf{SP} method uses 
finer grids: $N_v = 320$ for $\varepsilon = 0.2$, $N_v = 640$ for $\varepsilon = 0.1$, and $N_v = 1280$ for $\varepsilon = 0.05$.
We compute the relative $\ell^2$ errors computed by
\begin{equation*}
    \mathcal{E}^L_2 := \frac{\|f^L_{AE} - f^L_{\text{ref}}\|_{\ell^2_v}}{\|f^L_{\text{ref}}\|_{\ell^2_v}}, \quad \mathcal{E}^H_2 := \frac{\|f^H_{AE} - f^H_{\text{ref}}\|_{\ell^2_v}}{\|f^H_{\text{ref}}\|_{\ell^2_v}},
\end{equation*}
where $f^L_{\text{ref}}$, $f^H_{\text{ref}}$ denote the numerical solutions of $f^L$, $f^H$ solved by the \textbf{SP} method.
We set the time step size as $\Delta t = 0.1$. Table~\ref{table:compare_error} shows the results at the final time $t = 1$. We observe a decrease in relative errors as $\varepsilon$ becomes smaller, with a satisfactory level of accuracy of $\mathcal{O}(10^{-3})$ to $\mathcal{O}(10^{-5})$. This is consistent with the discussion in Section~\ref{subsec:AE_applicability} and meets our expectations. 
\begin{table}[!h]
    \centering
    \resizebox{0.7\textwidth}{!}{
    \begin{tabular}{c | cc | cc | cc}
    \toprule
     & \multicolumn{2}{c|}{$\epsilon = 0.2$} & \multicolumn{2}{c|}{$\epsilon = 0.1$} & \multicolumn{2}{c}{$\epsilon = 0.05$} \\
    \cmidrule(lr){2-3} \cmidrule(lr){4-5} \cmidrule(lr){6-7}
    $N_v$ & $\mathcal{E}^L_2$ & $\mathcal{E}^H_2$ & $\mathcal{E}^L_2$ & $\mathcal{E}^H_2$ & $\mathcal{E}^L_2$ & $\mathcal{E}^H_2$ \\
    \midrule
    40  & 8.81e-03 & 9.72e-03 & 7.80e-03 & 3.24e-03 & 7.57e-03 & 9.03e-04 \\
    80  & 4.13e-03 & 1.90e-03 & 3.87e-03 & 8.03e-04 & 3.81e-03 & 2.63e-04 \\
    160 & 2.02e-03 & 3.77e-03 & 1.93e-03 & 2.38e-04 & 1.93e-03 & 5.91e-05 \\
    320 & 1.01e-03 & 4.50e-03 & 9.79e-04 & 3.70e-04 & 9.82e-04 & 3.07e-05 \\
    \bottomrule
    \end{tabular}
    }
    \caption{\footnotesize
    Relative errors of the solution computed by the \textbf{AE} method compared with the \textbf{SP} method with different numbers of velocity points at time $t=1$.}
    \label{table:compare_error}
\end{table}

\begin{figure}
    \centering
        \includegraphics[width=0.8\textwidth]{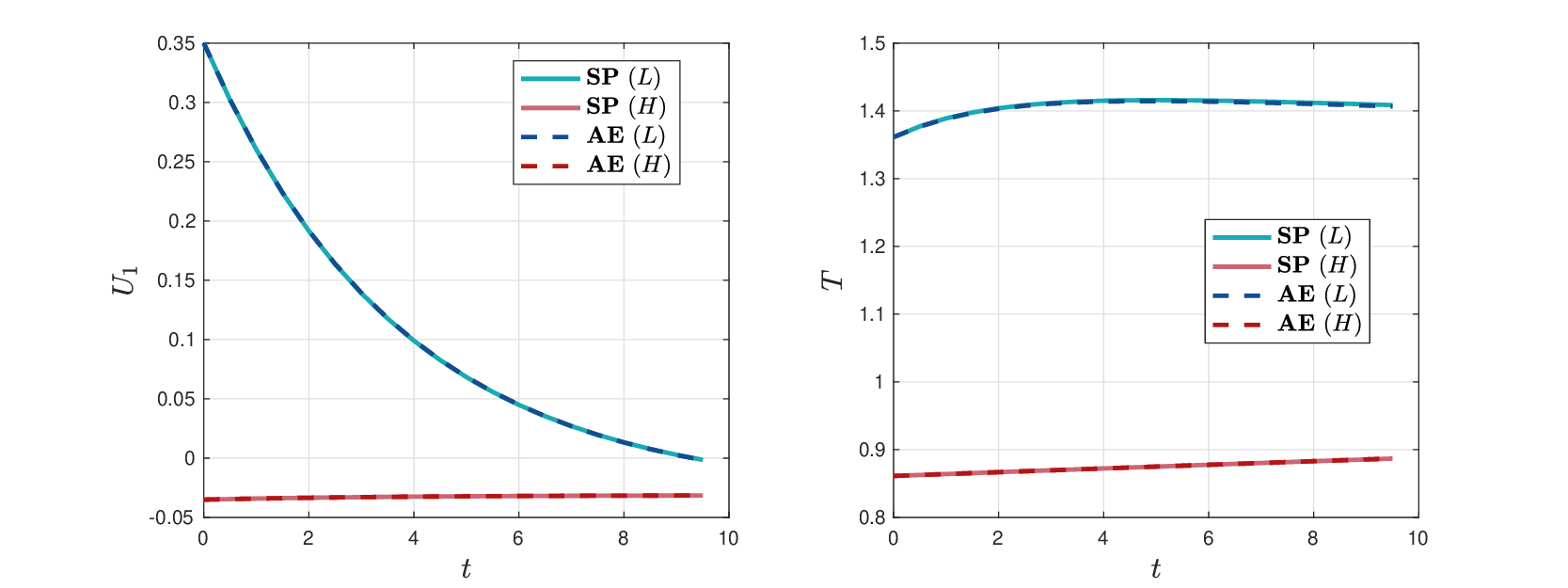}
    \caption{\footnotesize Evolution of the first direction of velocity ($U_1$) and temperatures ($T$) for $\varepsilon = 0.1$.}
    \label{figure:test1_evol}
\end{figure}

To compare the macroscopic quantities computed from the numerical solutions, we present the time evolution of the first direction of velocity vector and temperature for the case $\varepsilon = 0.1$, in Figure~\ref{figure:test1_evol}. The dashed lines represent the numerical solutions obtained by the \textbf{AE} method, while the solid lines represent the solution from the \textbf{SP} method. For the \textbf{AE} method, a velocity grid size of $N_v = 480$ is used, and the \textbf{SP} method uses $N_v = 800$. We set $\Delta t = 0.5$ in both methods.  

\subsection{Epochal relaxation}

In this test, we demonstrate the AP property of our scheme, namely to capture the asymptotic behaviors across multiple time scales in the disparate mass model~\eqref{equation:main_equation}. We validate this through the epochal relaxation phenomenon described in subsection~\ref{subsection:epochal_relaxation}. 

\vspace{10pt}

\begin{figure}
    \centering
    \begin{subfigure}{0.32\textwidth}
        \centering
        \includegraphics[width=1\textwidth]{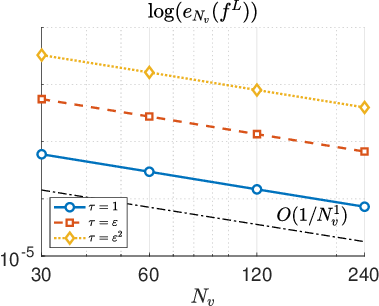}
    \end{subfigure}
    \hfill
    \begin{subfigure}{0.32\textwidth}
        \centering
        \includegraphics[width=1\textwidth]{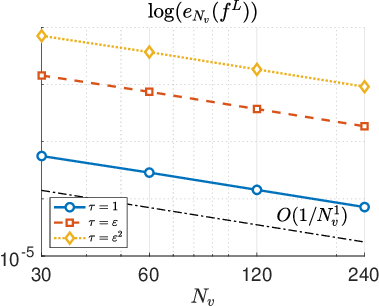}
    \end{subfigure}
    \hfill
    \begin{subfigure}{0.32\textwidth}
        \centering
        \includegraphics[width=1\textwidth]{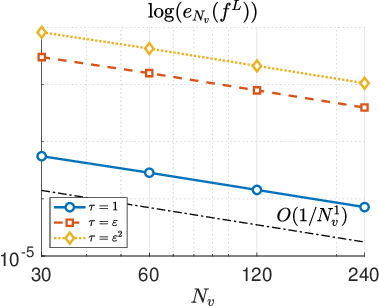}
    \end{subfigure}
    \vskip\baselineskip
    \begin{subfigure}{0.32\textwidth}
        \centering
        \includegraphics[width=1\textwidth]{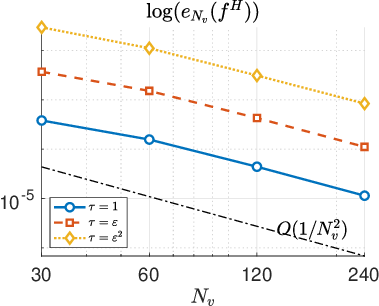}
        \caption*{\footnotesize $\varepsilon = 10^{-1}$}
    \end{subfigure}
    \hfill
    \begin{subfigure}{0.32\textwidth}
        \centering
        \includegraphics[width=1\textwidth]{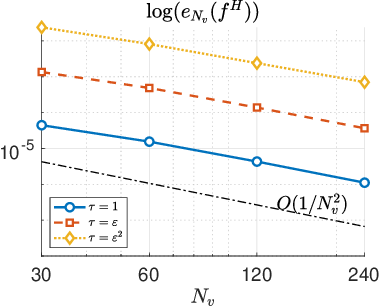}
        \caption*{\footnotesize $\varepsilon = 3 \times 10^{-2}$}
    \end{subfigure}
    \hfill
    \begin{subfigure}{0.32\textwidth}
        \centering
        \includegraphics[width=1\textwidth]{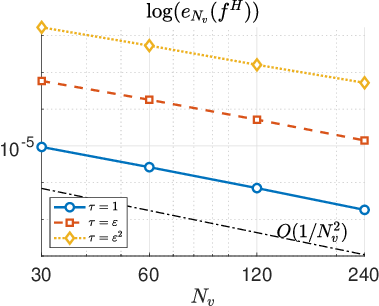}
        \caption*{\footnotesize $\varepsilon = 1 \times 10^{-2}$}
    \end{subfigure}
    \caption{\footnotesize Relative $\ell^2$ errors for the numerical solutions to the AP scheme~\eqref{equation:scheme} at the final time $t=0.1$.}
    \label{figure:accuracy}
\end{figure}

\noindent\textbf{Convergence in velocity of the AP scheme.} First we conduct the numerical computation for the AP scheme~\eqref{equation:scheme} for disparate mass gas mixtures, and adopt the \textbf{AE} method to approximate the inter-particle collision operators. We study gas mixture problems characterized by decreasing $\varepsilon$ values of $0.1$, 
$0.03$ and $0.01$ and at three different time scales. The convergence of numerical solutions in the velocity space will be investigated. The velocity points are chosen by $N_v=30$, $60$, $120$, $240$, and $480$ in all cases. We set the time step size as $\Delta t = 0.01$ and the final computation time $t = 0.1$. In Figure~\ref{figure:accuracy}, the relative $\ell^2$ errors of the solutions defined by~\eqref{equation:l2error} are shown. Similar to the result of Figure~\ref{figure:AE_SP_comparison}, we observe a second order accuracy for the numerical solution $f^H$ and first order accuracy for $f^L$, and they both enjoy a {\it uniform} accuracy in $\tau$.

Next, we investigate the asymptotic behaviors of the distribution functions and the macroscopic quantities at three time scales: $\tau = 1$, $\tau = \varepsilon$, and $\tau = \varepsilon^2$. 
We consider three disparate mass mixtures parameterized by $\varepsilon = 10^{-2}$, $10^{-3}$, and $10^{-4}$. 
For all simulations in the following parts, we use a time step size of $\Delta t = 0.1$ and a velocity grid size of $N_v = 200$.

\vspace{10pt}

\begin{figure}
    \centering
    \begin{subfigure}{0.4\textwidth}
        \centering
        \includegraphics[width=1\textwidth]{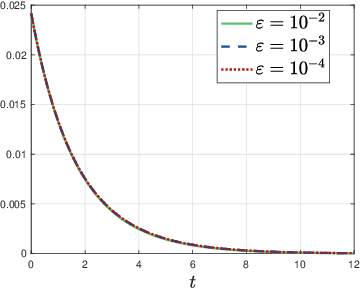}
        \label{figure:rel_entropyL_tau_1}
    \end{subfigure}
    \hfill
    \begin{subfigure}{0.4\textwidth}
        \centering
        \includegraphics[width=1\textwidth]{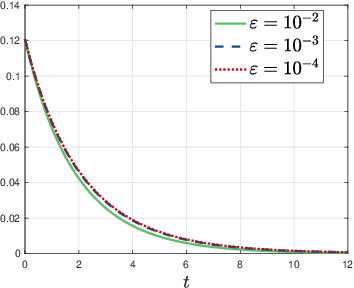}
        \label{figure:rel_entropyH_tau_eps}
    \end{subfigure}
    \caption{\footnotesize Separation of thermodynamic relaxation scales of the two species. Left: Time evolution of $H(f^L)$ for different $\varepsilon$ at the \textbf{fastest time scale ($\tau = 1$)}. Right: Time evolution of $H(f^H)$ for different $\varepsilon$ at the \textbf{intermediate time scale ($\tau = \varepsilon$)}. }
    \label{fig:entropy}
\end{figure}

\noindent\textbf{Maxwellization of distribution functions.} 
First, we observe a separation of thermodynamic relaxation scales for the light species and heavy species. More specifically, the Maxwellization of the light species happens at the fastest time scale ($\tau = 1$), while that of the heavy species happens at the intermediate time scale ($\tau = \varepsilon$). Figure~\ref{fig:entropy} shows the evolution of the relative entropies between the distribution functions and their local Maxwellians $M^L_0$ and $M^H$: 
\begin{equation*}
    H(f^L) = \int_{\mathbb{R}^2} f^L(v) \log \left(\frac{f^L(v)}{M^L_0(v)}\right)\mathrm{d} v, \quad
    H(f^H) = \int_{\mathbb{R}^2} f^H(v) \log \left(\frac{f^H(v)}{M^H(v)}\right)\mathrm{d} v, 
\end{equation*}
where $M^L_0$ denotes the centered local Maxwellian of $f^L$. 
According to~\eqref{equation:fL_leadingOrder} and~\eqref{equation:fH_leadingOrder}, the leading order behaviors of the two collision processes are the same for different $\varepsilon$. This is reflected by the similar relaxation behaviour for varying $\varepsilon$ in Figure~\ref{fig:entropy}. 
Moreover, at the fastest time scale, the light species particles go through elastic scattering against the heavy ones as if the latter were steady~\cite{Degond-Lucquin-Desreux}, and the zeroth order distribution of $f^L$ is a centered Maxwellian. 

\begin{figure}
    \centering
    \includegraphics[width=0.8\textwidth]{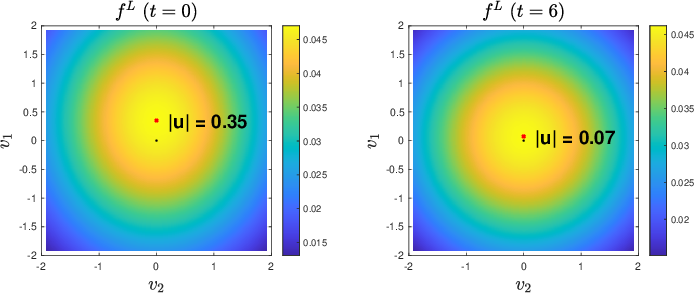}
    \label{figure:distribution_L_e0.01_tau_1}
    \caption{\footnotesize The light species distribution $f^L$ relaxes to a centered Maxwellian at the \textbf{fastest time scale ($\tau = 1$)}. Here $\varepsilon = 0.01$ and the red crosses denote the mean velocities. }
    \label{fig:dist_L}
\end{figure}

\begin{figure}
    \centering
    \includegraphics[width=0.8\textwidth]{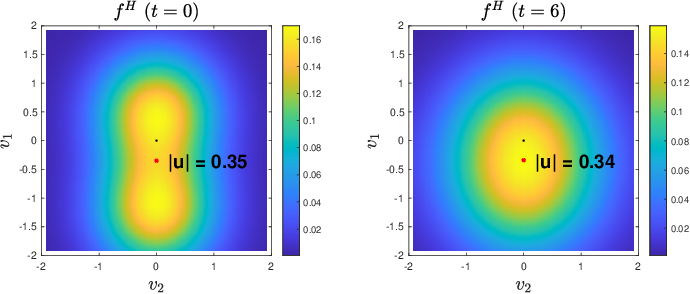}
    \label{figure:distribution_H_e0.01_tau_eps}
    \caption{\footnotesize The heavy species distribution $f^H$ relaxes to a non-centered Maxwellian at the \textbf{intermediate time scale ($\tau = \varepsilon$)}. Here $\varepsilon = 0.01$ and the red crosses denote the mean velocities. }
    \label{fig:dist_H}
\end{figure}

\begin{figure}
    \centering
    \begin{subfigure}{0.4\textwidth}
        \centering
        \includegraphics[width=1\textwidth]{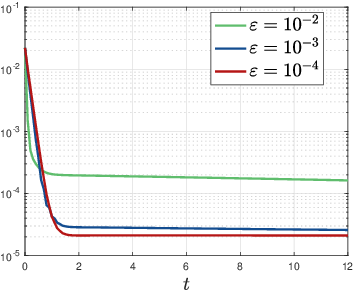}
        \label{figure:errorL_tau_eps}
    \end{subfigure}
    \hfill
    \begin{subfigure}{0.4\textwidth}
        \centering
        \includegraphics[width=1\textwidth]{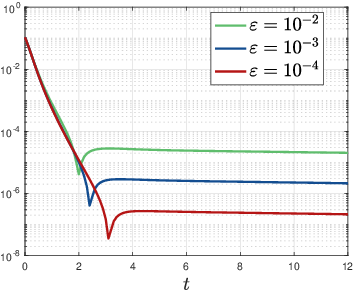}
        \label{figure:errorH_tau_eps2}
    \end{subfigure}
    \caption{\footnotesize Separation of hydrodynamic scales for the two species. Left: Time evolution of $\|f^L - M^L_0\|_{L^2}$ for different $\varepsilon$ at the \textbf{intermediate time scale ($\tau = \varepsilon$)}. Right: Time evolution of $\|f^H - M^H\|_{L^2}$ for different $\varepsilon$ at the \textbf{slowest time scale ($\tau = \varepsilon^2$)}. }
    \label{fig:APerror}
\end{figure}

Figure~\ref{fig:dist_L} and~\ref{fig:dist_H} show the snapshots of distribution solutions at time $t=0$ and $t=6$ for $f^L$ at the fastest time scale and $f^H$ at the intermediate time scale. From this figure, one can see that $f^L$ evolves towards a centered Maxwellian with $(0,0)$ velocities while $f^H$ evolves towards a non-centered one. Figure~\ref{fig:APerror} shows the time evolution of $\|f^L - M^L_0\|_{\ell^2}$ and $\|f^H - M^H\|_{\ell^2}$ at the intermediate time scale and the slowest time scale, which are the 
hydrodynamic scales of the light species and the heavy species, respectively. Under those time scales, the distribution functions relax quickly to their local Maxwellians. In addition, we find out that as $\varepsilon$ decreases, the errors when they reach a saturated level become smaller correspondingly. 
This phenomenon is consistent with the description in the previous work~\cite{Degond-Lucquin-Desreux} and our Table~\ref{tab:epochal}. 


\begin{figure}
    \centering
    \begin{subfigure}{0.32\textwidth}
        \centering
        \includegraphics[width=1\textwidth]{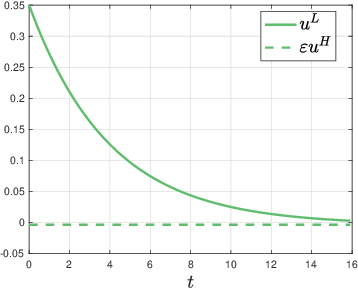}
        \caption*{\footnotesize $\varepsilon = 10^{-2}$}
        \label{figure:velocity_tau_1_e0.01}
    \end{subfigure}
    \hfill
    \begin{subfigure}{0.32\textwidth}
        \centering
        \includegraphics[width=1\textwidth]{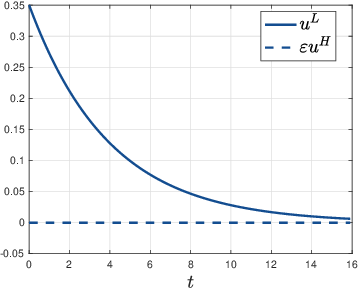}
        \caption*{\footnotesize $\varepsilon = 10^{-3}$}
        \label{figure:velocity_tau_1_e0.001}
    \end{subfigure}
    \hfill
    \begin{subfigure}{0.32\textwidth}
        \centering
        \includegraphics[width=1\textwidth]{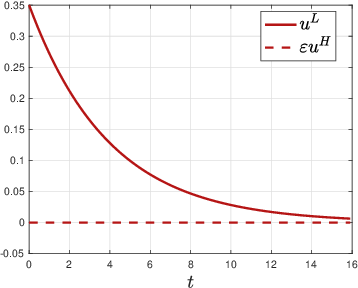}
        \caption*{\footnotesize $\varepsilon = 10^{-4}$}
        \label{figure:velocity_tau_1_e0.0001}
    \end{subfigure}
    \caption{\footnotesize Time evolution of $u^L$ and $\varepsilon u^H$ for different $\varepsilon$ at the \textbf{fastest time scale}; $u^L$ and $u^H$ are computed from numerical solutions $f^L$ and $f^H$ according to~\eqref{equation:moments}.}
    \label{fig:velocity}
\end{figure}

\begin{figure}
    \centering
    \begin{subfigure}{0.32\textwidth}
        \centering
        \includegraphics[width=1\textwidth]{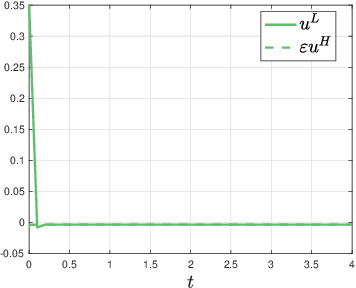}
        \caption*{\footnotesize $\varepsilon = 10^{-2}$}
        \label{figure:velocity_tau_1_e0.01}
    \end{subfigure}
    \hfill
    \begin{subfigure}{0.32\textwidth}
        \centering
        \includegraphics[width=1\textwidth]{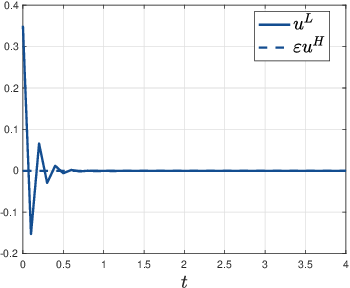}
        \caption*{\footnotesize $\varepsilon = 10^{-3}$}
        \label{figure:velocity_tau_1_e0.001}
    \end{subfigure}
    \hfill
    \begin{subfigure}{0.32\textwidth}
        \centering
        \includegraphics[width=1\textwidth]{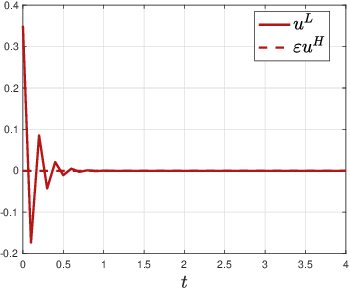}
        \caption*{\footnotesize $\varepsilon = 10^{-4}$}
        \label{figure:velocity_tau_1_e0.0001}
    \end{subfigure}
    \caption{\footnotesize Time evolution of $u^L$ and $\varepsilon u^H$ for different $\varepsilon$ at the \textbf{intermediate time scale}; $u^L$ and $u^H$ are computed from numerical solutions $f^L$ and $f^H$ according to~\eqref{equation:moments}.}
    \label{fig:velocity_tau_eps}
\end{figure}

\noindent\textbf{Relaxation of velocity and temperature.} According to \cite{Degond_bookchapter}, there is a separation of scales between the velocity relaxation and temperature relaxation of the species. In particular, the velocity relaxation happens at the fastest time scale, and the temperature relaxation happens at the slowest time scale. In Figure~\ref{fig:velocity} and~\ref{fig:velocity_tau_eps}, we show the evolution of mean velocities $u^L$ and $\varepsilon u^H$ for different $\varepsilon$ at the fastest time scale and the intermediate time scale. One can observe that the velocity relaxation happens at the fastest time scale. At the intermediate time scale, the difference between the velocities of the two species quickly become small. Meanwhile, their relaxation behaviours are similar for different regimes ($\varepsilon=10^{-2}$, $10^{-3}$ and $10^{-4}$). This result matches with our analysis in Table~\ref{tab:epochal}.

\begin{figure}
    \centering
    \begin{subfigure}{0.32\textwidth}
        \centering
        \includegraphics[width=1\textwidth]{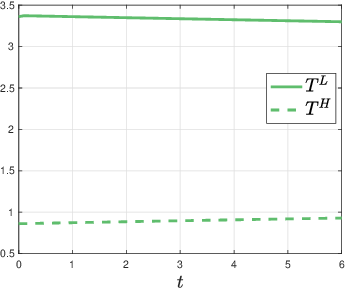}
        \label{figure:temperature_tau_eps2_e0.01}
        \caption*{\footnotesize $\varepsilon = 10^{-2}$}
    \end{subfigure}
    \hfill
    \begin{subfigure}{0.32\textwidth}
        \centering
        \includegraphics[width=\textwidth]{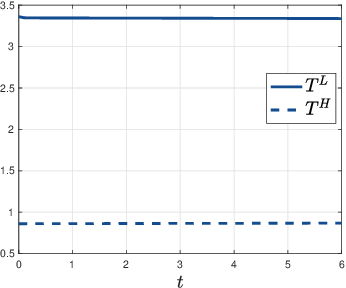}
        \label{figure:temperature_tau_eps2_e0.001}
        \caption*{\footnotesize $\varepsilon = 10^{-3}$}
    \end{subfigure}
    \hfill
    \begin{subfigure}{0.32\textwidth}
        \centering
        \includegraphics[width=\textwidth]{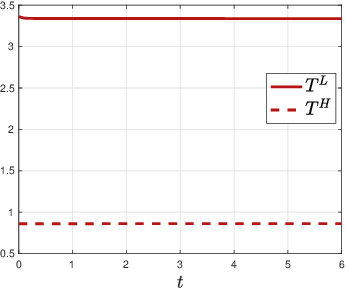}
        \label{figure:temperature_tau_eps2_e0.0001}
        \caption*{\footnotesize $\varepsilon = 10^{-4}$}
    \end{subfigure}
    \caption{\small Time evolution of temperatures for different $\varepsilon$ at the \textbf{intermediate time scale ($\tau = \varepsilon$)}. }
    \label{fig:temperature_tau_eps}
\end{figure}

\begin{figure}
    \centering
    \begin{subfigure}{0.32\textwidth}
        \centering
        \includegraphics[width=1\textwidth]{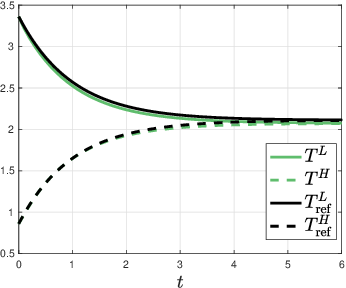}
        \label{figure:temperature_tau_eps2_e0.01}
        \caption*{\footnotesize $\varepsilon = 10^{-2}$}
    \end{subfigure}
    \hfill
    \begin{subfigure}{0.32\textwidth}
        \centering
        \includegraphics[width=1\textwidth]{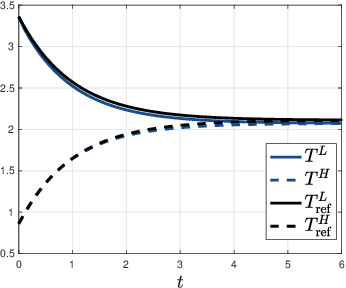}
        \label{figure:temperature_tau_eps2_e0.001}
        \caption*{\footnotesize $\varepsilon = 10^{-3}$}
    \end{subfigure}
    \hfill
    \begin{subfigure}{0.32\textwidth}
        \centering
        \includegraphics[width=1\textwidth]{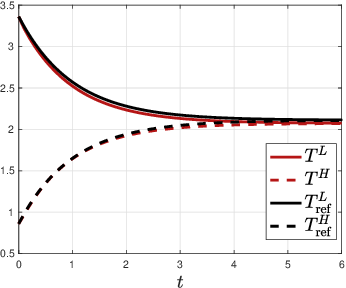}
        \label{figure:temperature_tau_eps2_e0.0001}
        \caption*{\footnotesize $\varepsilon = 10^{-4}$}
    \end{subfigure}
    \caption{\small Time evolution of temperatures for different $\varepsilon$ at the \textbf{slowest time scale ($\tau = \varepsilon^2$)}; $T^L$ and $T^H$ are computed from the numerical solutions $f^L$ and $f^H$ according to~\eqref{equation:moments}, and $T^L_\mathrm{ref}$ and $T^H_\mathrm{ref}$ are numerical solutions to~\eqref{equation:macro}.}
    \label{fig:temperature}
\end{figure}

Figure~\ref{fig:temperature_tau_eps}, on the other hand, shows the evolution of temperatures $T^L$ and $T^H$ computed by the distributions $f^L$ and $f^H$ at the intermediate time scale. As expected by our analysis from Table~\ref{tab:epochal}, the temperatures $T^L$ and $T^H$ remain constants at the leading order. Lastly, in Figure~\ref{fig:temperature}, we present the evolution of temperatures $T^L$ and $T^H$ computed by $f^L$ and $f^H$ at the slowest time scale, comparing with the solutions obtained from the macroscopic equations \eqref{equation:macro}. One can discover that the two macroscopic quantities match well, given the small $\varepsilon$, with the temperature of light and heavy species relaxing to each other as time is long enough.



\section{Conclusion}
\label{sec:conclusion}

In this paper, we have developed a new and efficient asymptotic-preserving (AP) scheme for the 
Boltzmann mixture model with disparate mass. A novel numerical method is constructed to compute the inter-particle collision operators on a {carefully constructed} polar-grid mesh, upon truncating their asymptotic expansions. The proposed numerical scheme remains {\it uniformly} accurate and computationally efficient for small mass ratio regimes. We carefully design the AP time-stepping in the moment updating step of the scheme, in order to ensure stability while capturing the epochal relaxation phenomenon without resolving the small scaling parameter. 
A series of numerical examples have demonstrated the strength and efficiency of our AP scheme across a range of  mass ratios. In future work, we will further address the accuracy and conservation preserving problem of the truncation method, in addition to study more general collision kernels and space-inhomogeneous problems.

\begin{appendices}

\section{The Boltzmann collision operators}
\label{appendix:collision_operator_formulae}

The intra-particle collision operators are the same as the single-species case~\cite{Cercignani, Chapman-Cowling}
{\small
\begin{equation}
    \label{equation:intra_collision_operators}
    \begin{aligned}
        &Q^{LL}(f^{L}, f^{L})(v^{L}) 
        = \int_{\mathbb{R}^{d_{v}}}\int_{\mathbb{S}^{d_{v}-1}}  B^{LL} (|v^{L} - v^{L}_*|, \sigma) 
        (f^{\prime L} f^{\prime L}_* - f^{L} f^{L}_*) \, 
        \mathrm{d}\sigma \, \mathrm{d}v^{L}_*, \\
        &Q^{HH}(f^{H}, f^{H})(v^{H}) 
        = \int_{\mathbb{R}^{d_{v}}}\int_{\mathbb{S}^{d_{v}-1}}  B^{HH} (|v^{H} - v^{H}_*|, \sigma) 
        (f^{\prime H} f^{\prime H}_* - f^{H} f^{H}_*) \, 
        \mathrm{d}\sigma \, \mathrm{d}v^{H}_*,
    \end{aligned}
\end{equation}
}
where $f^L = f^{L}(v^L)$, $f^H = f^{H}(v^H)$, $f^{L}_* = f^{L}(v^{L}_*)$, and $f^{H}_* = f^{H}(v^{H}_*)$.
The post-collision velocities can be parameterized through the collision transforms
{\small
\begin{align*}
    &v^{\prime L} = v^L
    + \frac{1}{2} |v^L - v^L_*| \sigma, \quad v^{\prime L}_* = v^L
    - \frac{1}{2} |v^L - v^L_*| \sigma \\
    &v^{\prime H} = v^H
    - \frac{1}{2} |v^H - v^H_*| \sigma, \quad v^{\prime H}_* = v^H
    + \frac{1}{2} |v^H - v^H_*| \sigma,
\end{align*}
}
The inter-particle collision operators are given by~\cite{Degond-Lucquin-Desreux}
{\small
\begin{equation}
    \begin{aligned}
        \label{equation:scaled_LH}
        Q^{LH}_{\varepsilon}(f^{L}, f^{H})(v^{L}) =& 
        \sqrt{1 + \varepsilon^2} \int_{\mathbb{R}^{d_{v}} \times \mathbb{S}^{d_{v}-1}} B^{LH}
        ( \tfrac{|g^{LH}|}{\sqrt{1 + \varepsilon^2}}, \sigma ) 
        ( f^{\prime L} f^{\prime H} - f^{L} f^{H} ) 
        \mathrm{d}\sigma \mathrm{d}v^{H}, 
    \end{aligned}
\end{equation}
\vspace{-10pt}
\begin{equation}
    \begin{aligned}
        \label{equation:scaled_HL}
        Q^{HL}_{\varepsilon}(f^{H}, f^{L})(v^{H}) =& 
        \tfrac{\sqrt{1 + \varepsilon^2}}{\varepsilon} 
        \int_{\mathbb{R}^{d_{v}} \times \mathbb{S}^{d_{v}-1}} B^{HL}
        (\tfrac{|g^{HL}|}{\sqrt{1 + \varepsilon^2}}, \sigma )
        ( f^{\prime H} f^{\prime L}- f^{H} f^{L} ) 
        \mathrm{d}\sigma \mathrm{d}v^{L}, 
    \end{aligned} 
\end{equation}
}
with collision rules
{\small
\begin{equation}
    \label{equation:collision_rules1}
        \begin{aligned}
            v^{\prime L} = v^{L} - \frac{1}{1 + \varepsilon^2} g^{LH} 
            + \frac{1}{1 + \varepsilon^2} |g^{LH}| \sigma, \quad
            \varepsilon v^{\prime H} = v^{L} - \frac{1}{1 + \varepsilon^2} g^{LH} 
            - \frac{\varepsilon^2}{1 + \varepsilon^2} |g^{LH}| \sigma. 
        \end{aligned}
\end{equation}
}
where $g^{LH} = v^L - \varepsilon v^{H} = -g^{HL}$.
$B^{LL}$, $B^{HH}$, $B^{LH}$, and $B^{HL}$ are the collision kernels. In the case of $d_v = 2$ and Maxwell molecules, the collision kernels are constants.

\section{Scaled spectral method}
  \label{appendix:SP}
  
  In this section, we provide a brief overview of the spectral method for evaluating the scaled operators $Q^{LH}_\varepsilon$ and $Q^{HL}_\varepsilon$. This approach is an adaptation of the fast spectral method developed by \cite{Jaiswal-Alexeenko-Hu}.
  
  \subsection{Re-scaling of operators}
  We first perform a set of variable changes
  {\small
  \begin{equation}
      \label{equation:change_variables}
      \tilde{v}^{H} = \varepsilon v^{H}, \quad f^{H}(v^{H}) = \varepsilon^{d_{v}} \tilde{f}^{H}(\tilde{v}^{H}), \quad Q^{HL}_{\varepsilon}(f^{H}, f^{L})(v^{H}) = \varepsilon^{d_{v}} \tilde{Q}^{HL}_{\varepsilon}(\tilde{v}^{H}),
  \end{equation}}
  \vspace{-10pt}
  to get
  \vspace{-10pt}
  {\small
  \begin{equation}
      \label{equation:QLH_unscaled}
      \begin{aligned}
          Q^{LH}_{\varepsilon}(v^{L}) =& \sqrt{1 + \varepsilon^2} \int_{\mathbb{R}^{d_{v}} \times \mathbb{S}^{d_{v}-1}} B^{LH}\left(\tfrac{|g^{LH}|}{\sqrt{1 + \varepsilon^2}}, \sigma \cdot \hat{g}^{LH}\right) \\
          & \cdot \left[ f^{L}(v^{\prime L}) \tilde{f}^{H}(\tilde{v}^{\prime H}) - f^{L}(v^{L}) \tilde{f}^{H}(\tilde{v}^{H}) \right] \mathrm{d}\sigma \mathrm{d}g^{LH},
      \end{aligned}
  \end{equation}}
  \vspace{-10pt}
    {\small \begin{equation}
      \label{equation:QHL_unscaled}
      \begin{aligned}
         \tilde{Q}^{HL}_{\varepsilon}(\tilde{v}^{H}) =& \tfrac{\sqrt{1 + \varepsilon^2}}{\varepsilon} \int_{\mathbb{R}^{d_{v}} \times \mathbb{S}^{d_{v}-1}} B^{HL}\left(\tfrac{|g^{HL}|}{\sqrt{1 + \varepsilon^2}}, \sigma \cdot \hat{g}^{HL}\right) \\
          & \cdot \left[ \tilde{f}^{H}(\tilde{v}^{\prime H}) f^{L}(v^{\prime L}) - \tilde{f}^{H}(\tilde{v}^{H}) f^{L}(v^{L}) \right] \mathrm{d}\sigma \mathrm{d}g^{HL},
      \end{aligned}
  \end{equation}
  }with collision rules~\eqref{equation:collision_rules1}.
  The integrations in \(\tilde{v}^{H}\) and \(v^{L}\) are transformed into integrations in \(g^{LH} = v^{L} - \tilde{v}^{H}\) and \(g^{HL} = \tilde{v}^{H} - v^{L}\), respectively.

  \subsection{Jaiswal-Alexeenko-Hu method~\cite{Jaiswal-Alexeenko-Hu}}
  
  With the new velocity variables $v^L$ and $\tilde{v}^H$, the Jaiswal-Alexeenko-Hu method can be applied to compute \eqref{equation:QLH_unscaled} and \eqref{equation:QHL_unscaled}. We set the velocity domain $\mathcal{D}_L = [-L,L]^2$ and periodize $f^L$ and $\tilde{f}^H$ to $\mathbb{R}^2$.
  Approximate the distribution functions by truncated Fourier series 
  \vspace{-5pt}
  {\small
  \begin{equation*}
      f^L_{N_v}(v^L) = \sum_{l = -\tfrac{N_v}{2}}^{\tfrac{N_v}{2} - 1} \hat{f}^L_l \, e^{\mathrm{i} \frac{\pi}{L} v^L \cdot l}, \quad 
      \tilde{f}^H_{N_v}(\tilde{v}^H) = \sum_{l = -\tfrac{N_v}{2}}^{\tfrac{N_v}{2} - 1} \hat{\tilde{f}}^H_l \, e^{\mathrm{i} \frac{\pi}{L} \tilde{v}^H \cdot l}.
  \end{equation*}}
  \vspace{-10pt}
  Substitute these into \eqref{equation:QLH_unscaled}, and project to the same Fourier space
  {\small
  \begin{equation*}
      \hat{Q}^{LH}_{\varepsilon, k} = 
      \frac{1}{(2L)^{d_{v}}}\sum_{\substack{l,m=-\tfrac{N_v}{2}\\l+m=k}}^{\tfrac{N_v}{2}-1}
      G^{LH}_{\varepsilon}(l,m)\hat{f}^{L}_{l}\hat{\tilde{f}}^{H}_{m}, \quad
      k \in \{-\tfrac{N_v}{2}, ... \tfrac{N_v}{2}-1\}.
  \end{equation*}}
  The kernel modes are given by $G^{LH}_{\varepsilon}(l,m) = G^{LH,+}_{\varepsilon}(l,m)-G^{LH,-}_{\varepsilon}(m)$,
  with 
  {\small
  \begin{align*}
   G^{LH,+}_{\varepsilon}(l,m) =& \sqrt{1+\varepsilon^2}\int_{\mathcal{B}_{R}\times\mathbb{S}^{d_{v}-1}}
      B^{LH}e^{-\tfrac{\mathrm{i}\pi}{L}\tfrac{(l+m)\cdot g^{LH}}{1+\varepsilon^2}
      + \tfrac{\mathrm{i}\pi}{L}\tfrac{|g^{LH}|\sigma\cdot(l - \varepsilon^2 m)}{1+\varepsilon^2}}
      \mathrm{d}\sigma\mathrm{d}g^{LH},\\
      G^{LH,-}_{\varepsilon}(m) =& \sqrt{1+\varepsilon^2}\int_{\mathcal{B}_{R}\times\mathbb{S}^{d_{v}-1}}
      B^{LH}e^{-\tfrac{\mathrm{i}\pi}{L}m \cdot g^{LH}}
      \mathrm{d}\sigma\mathrm{d}g^{LH}.
  \end{align*}}The integration in $g^{LH}$ is dealt with in spherical coordinates.
  In the case where $d_v=2$ and Maxwell molecules, we can derive
  $\hat{Q}^{LH}_{\varepsilon, k} = \hat{Q}^{LH,+}_{\varepsilon, k} - \hat{Q}^{LH,-}_{\varepsilon, k}$ and $\hat{\tilde{Q}}^{HL}_{\varepsilon, k} = \hat{\tilde{Q}}^{HL,+}_{\varepsilon, k} - \hat{\tilde{Q}}^{HL,-}_{\varepsilon, k}$, where
  {\small
  \begin{align*}
      \hat{Q}^{LH,+}_{\varepsilon, k} =& 2\pi B^{LH} \sqrt{1 + \varepsilon^2} \sum_{\rho, \theta}\sum_{l + m = k}  w_{\rho} w_{\theta} \rho \mathcal{J}_{0}\left(\tfrac{\pi}{L} \rho \tfrac{|k|}{1 + \varepsilon^2}\right) \hat{f}^{L}_{l} \, \hat{\tilde{f}}^{H}_{m} \, e^{\mathrm{i} \frac{\pi}{L} \rho \frac{|l| \cos \theta - \varepsilon^2 |m| \cos \theta}{1 + \varepsilon^2}}, \\
      \hat{\tilde{Q}}^{HL,+}_{\varepsilon, k} =& 2\pi B^{HL} \frac{\sqrt{1 + \varepsilon^2}}{\varepsilon} \sum_{\rho, \theta}\sum_{l + m = k}  \rho \mathcal{J}_{0}\left(\tfrac{\pi}{L} \rho \tfrac{\varepsilon^2 |k|}{1 + \varepsilon^2}\right)  \hat{\tilde{f}}^{H}_{l} \, \hat{f}^{L}_{m} \, e^{\mathrm{i} \frac{\pi}{L} \rho \frac{\varepsilon^2 |l| \cos \theta - |m| \cos \theta}{1 + \varepsilon^2}}, \\
      \hat{Q}^{LH,-}_{\varepsilon, k} =& \sum_{l + m = k}^{\tfrac{N_v}{2} - 1} \hat{f}^{L}_{l} \left[G^{LH,-}_{\varepsilon}(m) \hat{\tilde{f}}^{H}_{m}\right], \; \hat{\tilde{Q}}^{HL,-}_{\varepsilon, k} = \sum_{l + m = k}^{\tfrac{N_v}{2} - 1} \hat{\tilde{f}}^{H}_{l} \left[G^{HL,-}_{\varepsilon}(m) \hat{f}^{L}_{m}\right]. 
  \end{align*}}
  The expressions of $G^{LH,-}_{\varepsilon}(m)$ and $G^{HL,-}_{\varepsilon}(m)$ are given by
  {\small
  \begin{align*}
      G^{LH,-}_{\varepsilon}(m) =& 
      \left\{
          \begin{array}{l}
          4\pi^2 R^2 B^{LH}\sqrt{1+\varepsilon^2}
          \mathcal{J}_{1}(\tfrac{R\pi}{L}|m|)/\tfrac{R\pi}{L}|m|,\quad m\neq 0\\
          2\pi^2 R^2 B^{LH}\sqrt{1+\varepsilon^2},\quad m=0.
      \end{array} 
      \right. \\
      G^{HL,-}_{\varepsilon}(m) =&  
      \left\{
          \begin{array}{l}
          4\pi^2 R^2 B^{HL} \tfrac{\sqrt{1 + \varepsilon^2}}{\varepsilon} \mathcal{J}_{1}\left(\tfrac{R \pi}{L} |m|\right) / \left(\tfrac{R \pi}{L} |m|\right), \quad m \neq 0 \\
          2\pi^2 R^2 B^{HL} \tfrac{\sqrt{1 + \varepsilon^2}}{\varepsilon}, \quad m = 0,
          \end{array}
      \right.
  \end{align*}}
  where $\mathcal{J}_0(r)$, $\mathcal{J}_1(r)$ are Bessel functions.
  
  \subsection{Fourier transforms for scaled functions}
    To numerically handle the re-scaling \( f^H(v^H) = \varepsilon^{d_v} \tilde{f}^H(\tilde{v}^H) \)
    we leverage the dilation property of Fourier transforms to handle the numerical rescaling. Specifically, by assuming that \(\tilde{f}^{H}(\tilde{v}^{H})\) has support \(\mathcal{B}_{\varepsilon S}\), we have by~\eqref{equation:change_variables}, 
    {\small
  \begin{equation}
      \label{equation:FT}
        \hat{\tilde{f}}^{H}_{m} = \frac{1}{(2L)^{d_{v}}} \int_{\mathcal{D}_{L}} \tilde{f}^{H}(\tilde{v}^{H}) e^{-\mathrm{i} \frac{\pi}{L} m \cdot \tilde{v}^{H}} \mathrm{d}\tilde{v}^{H}
          = \frac{1}{(2L)^{d_{v}}} \int_{\mathcal{D}_{L}} f^{H}(v^{H}) e^{-\mathrm{i} \frac{\pi}{L} \varepsilon m \cdot v^{H}} \mathrm{d}v^{H}. 
  \end{equation}
}We then obtain the coefficients $\{\hat{\tilde{f}}^H_m\}$ directly from an integration of $f^H$, without resorting to the numerical representation of $\tilde{f}^H$. Similarly, we can get $Q^{HL}$ directly from $\{\hat{\tilde{Q}}^{HL}_{\varepsilon, k}\}$ by
{\small $ Q^{HL}_{\varepsilon}(f^{H}, f^{L})(v^{H}) = \frac{\varepsilon^2}{(2L)^2} \sum_{\substack{k = -\tfrac{N_v}{2}}}^{\tfrac{N_v}{2} - 1} \hat{\tilde{Q}}^{HL}_{\varepsilon, k} e^{\mathrm{i} \frac{\pi}{L} \varepsilon k \cdot v^{H}}. $ }

\subsection{Time complexity of the \textbf{SP} method}

Based on the above algorithm, we can derive the time complexity of the \textbf{SP} method.
\begin{proposition}
\label{Prop:SP}
The \textbf{SP} method has a time complexity of $\mathcal{O}\left( (N_v(\varepsilon))^{2d_v} \right)$, with $N_v(\varepsilon) \propto 1/\varepsilon$. 

\begin{proof}
    In the \textbf{SP} method, we assume 
$ \text{Supp}(\tilde{f}^{H})
\subset \mathcal{B}_{\varepsilon S}$. Here, $\mathrm{Supp}(\tilde{f}^H)$ denotes the support of $\tilde{f}^H$, and $\mathcal{B}_{S}$ is the ball in $\mathbb{R}^{d_v}$ centered at the origin with radius $S$. By the Heisenberg's inequality \cite{folland_fourier}, one has 
    $|\mathrm{Supp}(\tilde{f}^H)| \cdot |\mathrm{Supp}(\hat{\tilde{f}}^H_m)| \gtrsim 1 $, 
    where $|\cdot|$ denotes the measure of a set, $\tilde{f}^H$ and $\hat{\tilde{f}}^H_m$ are defined in \eqref{equation:FT}. Thus  $|\mathrm{Supp}(\hat{\tilde{f}}^H_m)| \approx 1/\varepsilon. $
    To accurately capture this change in computational domain from $\tilde{f}^H$ to $\hat{\tilde{f}}^H_m$, for the velocity discretization $N_v(\varepsilon) \propto 1/\varepsilon$ points are needed. 
\end{proof}
\end{proposition}

  \section{Derivation of asymptotic operators}
  \label{appendix:derivation}
  The derivation is based on the Taylor expansion of the integrand, following the collision rule given by~\eqref{equation:collision_rules1}
  {\small  
  \begin{equation*}
      f^L(v^{\prime L}) f^H(v^{\prime H}) - f^L(v^L) f^H(v^H) = I_0 + \varepsilon I_1 + \varepsilon^2 I_2 + \mathcal{O}(\varepsilon^3), 
  \end{equation*}
  }
  where
  {\footnotesize
  \begin{align*}
      I_0 =& \left( f^L(|v^L| \sigma) - f^L(v^L) \right) f^H(v^H), \\
      I_1 =&  f^L(|v^L| \sigma) (v^L - |v^L| \sigma) \cdot \nabla_{v^H} f^H(v^H) + f^H(v^H) \big(v^H - \frac{(v^L \cdot v^H)}{|v^L|} \sigma\big) 
      \cdot \nabla_{v^L} f^L(|v^L| \sigma), \\
      I_2 =& \frac{1}{2} f^L(|v^L| \sigma) (v^L - |v^L| \sigma)^{\otimes 2} : \nabla_{v^H}^2 f^H(v^H) + f^L(|v^L| \sigma) \big(-v^H + \frac{(v^L \cdot v^H)}{|v^L|} \sigma\big) 
          \cdot \nabla_{v^H} f^H(v^H) \\
          &+ \frac{1}{2} f^H(v^H) \big(v^H - \frac{(v^L \cdot v^H)}{|v^L|} \sigma\big)^{\otimes 2}: \nabla_{v^L}^2 f^L(|v^L| \sigma) \\
      &+ \big(v^H - \frac{(v^L \cdot v^H)}{|v^L|} \sigma \big)
      \cdot \nabla_{v^L} f^L(|v^L| \sigma) (v^L - |v^L| \sigma)
      \cdot \nabla_{v^H} f^H(v^H) \\
      &+ f^H(v^H) \big(v^L + \frac{1}{2} |v^L| \sigma 
      \big(\frac{|v^H|^2}{|v^L|^2} - \frac{(v^L \cdot v^H)^2}{|v^L|^4}\big) - |v^L| \sigma\big)
      \cdot \nabla_{v^L} f^L(|v^L| \sigma), 
  \end{align*}
  }
  To deal with the singularity at the origin, we let $v^L=0$ in~\eqref{equation:collision_rules1}. We get $I_0 = 0$ and
  {\small
  \begin{align*}
        I_1 =& (v^H + |v^H|\sigma)\cdot\nabla_{v^L} f^L(0) f^H(v^H), \;
      I_2 = f^L(0) (-v^H - |v^H| \sigma) \cdot \nabla_{v^H} f^H (v^H). 
  \end{align*}
  }
\end{appendices}

\bibliographystyle{siamplain}
\bibliography{Ref.bib}

\end{document}